\documentclass[12pt]{article}

\usepackage
[colorlinks=true, pdfstartview=FitV, linkcolor=blue, citecolor=blue, urlcolor=blue]
{hyperref}

\usepackage{amssymb,amsmath, amscd,xypic, amsxtra}
\usepackage[T1]{fontenc}
\usepackage{times, verbatim}
\usepackage{graphicx}
\input xy
\xyoption{all}
\usepackage{tikz}
\usepackage{tikzsymbols}
\usepackage{mathtools}
\usepackage{tikz-3dplot}
\usepackage{pgfplots}
\usepackage{wrapfig}
\usepackage{caption}
\usepackage{skull}
\usepgfplotslibrary{polar}
\usepackage{xcolor}

\usepackage{mathpazo} % math & rm
\usepackage[scaled]{helvet} % ss
\usepackage{courier} % tt
\normalfont
\usepackage[T1]{fontenc}
\usepackage{esint}
\usepackage{verbatim}
\usetikzlibrary{arrows.meta}
\usetikzlibrary{arrows,decorations.markings}
\tikzset{>=latex}

\newtheorem{thm}{Theorem}

\newtheorem{lemma}[thm]{Lemma}
\newtheorem{cor}[thm]{Corollary}

\newcommand\proof
{\par\noindent{\bf Proof:\ } }
\newcommand\qed{\hfill$\blacksquare$}

\newcommand\cox{\mathsf{cox}}

\newcommand\blt{\bullet}

\newcommand{\wh}{\widehat}
        
\DeclareMathOperator{\Ad}{Ad}

\DeclareMathOperator{\Aut}{Aut}
\DeclareMathOperator{\ad}{ad}

\DeclareMathOperator{\Irr}{Irr}

\DeclareMathOperator{\tr}{tr}

\DeclareMathOperator{\Out}{Out}
\DeclareMathOperator{\rank}{rank}

\newcommand{\sw}{\mathsf{sw}}

\newcommand{\br}{\mathbb{R}}
\newcommand{\bc}{\mathbb{C}}

\newcommand{\bq}{\mathbb{Q}}
\newcommand{\bz}{\mathbb{Z}}

\newcommand{\al}{\alpha}
\newcommand{\be}{\beta}
\newcommand{\ga}{\gamma}

\newcommand{\Ga}{\Gamma}
\newcommand{\lam}{\lambda}

\newcommand{\Om}{\Omega}
\newcommand{\ep}{\epsilon}

\newcommand{\vp}{\varphi}

\newcommand{\scd}{\mathcal{D}}

\newcommand{\sci}{\mathcal{I}}

\newcommand{\scw}{\mathcal{W}}
\newcommand{\scx}{\mathcal{X}}
\newcommand{\scy}{\mathcal{Y}}
\newcommand{\scz}{\mathcal{Z}}

\newcommand{\la}{\langle}
\newcommand{\ra}{\rangle}

\DeclareMathOperator{\SO}{SO}

\DeclareMathOperator{\sgn}{sgn}

\newcommand{\bG}{\mathbf{G}}

\newcommand{\fg}{\mathfrak{g}}

\newcommand{\fsl}{\mathfrak{sl}}
\newcommand{\fso}{\mathfrak{so}}
\newcommand{\fsp}{\mathfrak{sp}}

\newcommand{\ft}{\mathfrak{t}}

\newcommand{\EVI}[7]{
{\scriptsize\begin{matrix}\!\!\!\!#1\!\!\!\!&\!\!\!\!#2\!\!\!\!&\!\!\!#3\!\!\!&\!\!\!\!#4\!\!\!\!
&\!\!\!\!#5
\!\!\!\!\\
&&\!\!\!#6\!\!\!&&\\&&\!\!\!#7\! \!\!&&\end{matrix}} }

\newcommand{\EVII}[8]{
{\scriptsize\begin{matrix}
 #1 & \!\!#2 \!& \!#3 \!\!& \!\!#4 \!\!& \!\!#5 \!& \!#6
 \!& \!#7\\
&&& #8 &&&\end{matrix}} }

\newcommand{\E}[9]{
{\scriptsize\begin{matrix}
\!\!\!#1\!\!\! \!\!&\!\!\!\!\!#2\!\!\!\!\!&\!\!\!\!\!#3\!\!\!\!\!&\!\!\!\!\!#4\!\!\!\!\!&\!\!\!\!\!#5\!\!\!\!\!&\!\!\!\!\!#6\!\!\!\!\!&\!\!\!\!\!#7
\!\!\!\!\!&\!\!\!\!\!#8\!\!\!\\
&&&&&\!#9\!\!\!&&\end{matrix}} }

\newcommand{\twoAtwo}{
\SelectTips{eu}{}
\xymatrix@M=0pt{
\underset{\hphantom{1}}{\overset{1}{\circ}}\ar@{=>}[r]\ar@<0.65ex>@{-}[r]+<-1.2ex,0ex>\ar@<-0.65ex>@{-}[r]+<-1.2ex,0ex>&
\underset{\hphantom{2}}{\overset{2}\circ}}
}

\newcommand{\twoAtwoo}{
\SelectTips{eu}{}
\xymatrix@M=0pt{
{1\ }\ar@{=>}[r]\ar@<0.65ex>@{-}[r]+<-1.2ex,0ex>\ar@<-0.65ex>@{-}[r]+<-1.2ex,0ex>&{\ 1}}
}
\newcommand{\twoAtwooo}{
\SelectTips{eu}{}
\xymatrix@M=0pt{
{1\ }\ar@{=>}[r]\ar@<0.65ex>@{-}[r]+<-1.2ex,0ex>\ar@<-0.65ex>@{-}[r]+<-1.2ex,0ex>&{\ 0}}
}

\newcommand{\outEVI}[5]{#1 #2 #3\!\Leftarrow\!#4 #5}
\newcommand{\FIV}[5]{#1 #2 #3\!\Rightarrow\!#4 #5}

\title {Thomae's function on a Lie group}

\author{Mark Reeder\\
Department of Mathematics, Boston College \\ 
Chestnut Hill, MA  02467\\
\texttt{reederma@bc.edu}
}

\date{\today}
\begin{document}
\maketitle
\begin{abstract}
Let $\fg$ be a simple complex Lie algebra of finite dimension. 
This paper gives an inequality relating the order of an automorphism of $\fg$ to the dimension of its fixed-point subalgebra, and characterizes those automorphisms of $\fg$ for which equality occurs. This is amounts to an inequality/equality for Thomae's function on $\Aut(\fg)$. The result has applications to characters of zero weight spaces, graded Lie algebras, and inequalities for adjoint Swan conductors.

\end{abstract}
\tableofcontents

%\date{\today}

\section{Introduction}

Thomae's function $\tau:\br\to\br$  is discontinuous  precisely on the rational numbers. It is traditionally defined as $\tau(x)=1/m$ if $x=n/m$ is rational in lowest terms with $m>0$, and $\tau(x)=0$ if $x$ is irrational. 
So $\tau(n)=1$ for every integer $n$, and  
on each open interval $(n,n+1)$ the maximum value of  $\tau$ is $1/2$, taken just at the midpoint of the interval.
More succinctly, $\tau(x)$ is the reciprocal of the order of $x$ in the group $\br/\bz$,
with the convention that $1/\infty=0$.  

Every  group $G$ has an analogous function $\tau_G:G\to\bq$, whose value at $g\in G$ equal to the reciprocal of the order of $g$.  

Consider the group $G=\SO_3$ of rotations about a fixed point $O$ in three-dimensional Euclidean space. Here $\tau_G(g)=1/m$ if $g$ rotates by a rational multiple $n/m$ (in lowest terms) of a full circle and $\tau_G(g)=0$ otherwise. So $\tau_G(g)=1$ if $g$ is the identity rotation and elsewhere $\tau_G$ has maximum value $1/2$ taken just on the conjugacy class of half-turns. Since every element of $G$ is conjugate to a rotation about a fixed axis through $O$, this example is essentially the same as Thomae's original one, but now we observe that $1/2=1/h$, where $h$ is the Coxeter number of $G$.

Suppose $G$ is either a compact Lie group or a complex algebraic group. For such groups the function $\tau_G$ is discontinuous precisely on the set of torsion elements in $G$. 
The proof is the same as for $\tau=\tau_{\br/\bz}$, using the facts 1) that torsion elements can be approximated by elements of infinite order, 2) for every $\ep>0$ there are only finitely many conjugacy classes in $G$ whose elements have order $\leq 1/\ep$, and 3) the conjugacy class of any torsion element is closed in $G$.

If $G$ is connected and simple as an abstract group then on the regular elements of $G$ we have $\tau_G(g)\leq 1/h$, where $h$ is the Coxeter number of $G$. Equality holds on just the conjugacy class of {\it principal elements}. These are the analogues of the half-turns in $\SO_3$ and were studied by Kostant in \cite{kostant:betti}. 

The aim of this paper is to extend this inequality/equality for Thomae's function to singular elements in the group $G=\Aut(\fg)$ of automorphisms of a simple complex Lie algebra $\fg$ of finite dimension. We also indicate some applications of the result. 

 We will measure the singularity of an element $\theta\in G$ by the dimension of the fixed-point subalgebra $\fg^\theta$. We  will give an upper bound for $\tau_{G}(\theta)$ in terms of $\dim\fg^\theta$, along with precise conditions for equality. 

To explain these conditions we need some preparation.  We say that an element $\theta\in G$ is {\bf ell-reg} if $\theta$ normalizes a Cartan subalgebra $\ft\subset \fg$ such that (i) $\ft^\theta=0$ and 
(ii) the cyclic group generated by $ \theta$ permutes the roots of $\ft$ in $\fg$ freely. There are only finitely many ell-reg classes in $\Aut(\fg)$. Their classification was given in \cite{rlyg} and is recalled in the appendix to this paper. \footnote{Such automorphisms are called {\it $\bz$-regular } in \cite{rlyg}, in deference to \cite{springer:regular}.  In general, ell-reg elements  are not regular elements of $G$. The point of ``ell-reg", besides brevity, is to avoid conflict between these two meanings of the word ``regular".}

For ell-reg automorphisms it is known that 
the automorphism of $\ft$ given by $\theta\vert_\ft$ (as in (i) and (ii)) has the same order as $\theta$. 
It follows  that if $\theta\in G$ is ell-reg then 
\begin{equation}\label{ell-reg}
\tau_{G}(\theta)= \frac{\dim\fg^\theta}{\dim(\fg/\ft)},
\end{equation}
where $\ft$ is any Cartan subalgebra of $\fg$. 

Fix a connected component $\Ga$ of $G$ and let $e\in\{1,2,3\}$ be the order of $\Ga$ in the  group $\Out(\fg)$ of connected components of $G$. 
If $\theta\in \Ga$, the rank of $\fg^\theta$ depends only on $e$; we write 
\[n_e=\rank(\fg^\theta).\]
In $\Ga$ there is a unique conjugacy class $P_\Ga$ of 
elements $\theta$ of minimal order for which $\fg^\theta$ is a Cartan subalgebra of $\fg^\theta$. This order, denoted $h_e$, is the {\it twisted Coxeter number} of the coset $\Ga$ \cite{reeder:torsion}. 
The elements of $P_\Ga$ are ell-reg and it is known that 
\begin{equation}\label{he}
\frac{1}{h_e}=\frac{n_e}{\dim(\fg/\ft)},\quad\text{if}\quad \theta \in P_\Ga.
\end{equation}
It follows that if $\theta\in \Ga$ has order $m>h_e$, then 
\begin{equation}\label{>he}
\tau_{G}=\frac{1}{m}<\frac{\dim\fg^\theta}{\dim(\fg/\ft)},
\end{equation}
Where $\tau_{G}$ is Thomae's function for the group $G=\Aut(\fg)$. 
In this paper we extend \eqref{>he} to all $\theta\in\Aut(\fg)$ as follows. 

\begin{thm}\label{mainthm} Let $\fg$ be a simple complex Lie algebra of finite dimension and let $\tau_{G}$ be Thomae's function for the group $G=\Aut(\fg)$. 
Then for all $\theta\in G$ we have 
\begin{equation}\label{ineq}
\tau_{G}(\theta)\leq \frac{\dim\fg^\theta}{\dim(\fg/\ft)}.
\end{equation}
Equality holds in \eqref{ineq} if and only if $\theta$ is ell-reg. 
\end{thm}

From \eqref{he} we have equality in \eqref{ineq} if $\theta\in P_\Ga$. Also \eqref{ineq} holds trivially, and is a strict inequality, if the order of $\theta$ is larger than $h_e$, by \eqref{>he}. 
Therefore, the content of Theorem 1 is (i) the inequality \eqref{ineq} for all $\theta\in G$ 
whose order $m$ lies in the range $1<m<h_e$, and (ii)
the assertion  that only ell-reg elements attain equality.

The proof of Theorem 1 consists of computations with Kac diagrams. 
It is given in section \ref{proof}.

\section{Applications}

First we give some applications of Theorem 1 and connections  to other results.  

\subsection{Characters of zero-weight spaces}\label{zeroweight}

This section describes the role played by Theorem 1 in the computation of characters of zero weight spaces \cite{reeder:weyl}.
\footnote{An earlier version of this paper was an appendix to an earlier version of \cite{reeder:weyl}.} 

Let $G$ be a connected and simply connected complex Lie group. Let $T$ be a maximal torus in $G$, with normalizer $N$ and Weyl group $W=N/T$. 
In every finite-dimensional irreducible representation  $V$ of $G$, the zero weight space $V^T$
affords a representation of $W$. The problem is to compute the character of $V^T$ afforded by $W$, as a function of the highest weight of $V$. 

For example, in \cite{kostant:eta}, Kostant used his results on principle elements to calculate the  trace $\tr(\cox,V^T)$ of a Coxeter element $\cox\in W$. 
He showed that   $\tr(\cox,V^T)$ is $0$ or $\pm 1$ and 
gave an explicit formula for this trace in terms of 
the highest weight of $V$. 

In \cite{prasad:cox}, Kostant's proof was reformulated in terms of the dual group $\wh G$ of $G$. Since $G$ is simply-connected, $\wh G$ is the group of inner automorphisms of the Lie algebra $\wh\fg$ whose root system is dual to that of $\fg$. In \cite{reeder:weyl}, Theorem 1 is used in $\wh G$ to compute traces of other Weyl group elements on $V^T$.

We call an element $w\in W$ {\it ell-reg} if (i) $\ft^w=0$ and (ii) the group $\la w\ra$ generated by $w$ acts freely on the roots of $\ft$ in $\fg$. Equivalently, $w$ is ell-reg if $\Ad(n)$ is an ell-reg automorphism of $\fg$ for one (every) lift of $w$ in $N$. The classification of ell-reg elements of $W$ is given by the untwisted diagrams in the Appendix \ref{appendix}.

Let  $w\in W$ be ell-reg of order $m$ and let $V=V_\lam$ have highest weight $\lam$, with respect to a choice of positive roots $R^+$ for $T$ in $G$.
Let $P$ and $Q$ be the weight and root-lattices of $T$, so that $P$ and $Q$ are also the co-weight and co-root lattices of the dual maximal torus $\wh T$ of $\wh G$. Let $\rho\in P$ be the half-sum of the roots in $R^+$. 
Let $\theta\in\wh T$ be the value at $e^{2\pi i/m}$ of the one-parameter subgroup $\lam+\rho$.

From Theorem 1 we get an inequality of centralizers 
\begin{equation}\label{d(t)}
\dim C_G(t)\leq \dim C_{\wh G}(\theta),
\end{equation}
with equality if and only if $(\lam+\rho)+mQ$ is conjugate to $\rho+mQ$ under the natural $W$-action on $P/mQ$. 
%This implies that the sum $\tr(w,V)$ of harmonic polynomials arising from $C_G(t)$ is actually single monomial arising from $C_{\wh G}(\theta)$. 
From this inequality and the theory of $W$-harmonic polynomials we deduce that
$\tr(w,V_{\lam}^T)=0$ unless there exists $v\in W$ such that $v(\lam+\rho)\in \rho+mQ$, in which case
\[
\tr(w,V_{\lam}^T)=\sgn(v)\prod_{\check\al\in \check R_m^+}
\frac{\la v(\lam+\rho),\check\al\ra}{\la \rho,\check\al\ra},
\]
where the product is over the positive co-roots $\check\al$ of $G$ for which $\la \rho,\check\al\ra\in m\bz$. 
If $m=h$ is the Coxeter number then $\check R^+_m=\varnothing$ and we recover Kostant's result for $\tr(\cox, V_\lam^T)$.

\subsection{Graded Lie algebras} 
Let $\theta\in\Aut(\fg)$ have order $m$, and let $\zeta=e^{2\pi i/m}$. 
Then $\theta$ determines a $\bz/m\bz$ grading 
\begin{equation}\label{grading}
\fg=\bigoplus_{k\in \bz/m\bz} \fg_k\ ,
\end{equation}
where $\fg_k=\{x\in \fg:\ \theta(x)=\zeta^k x\}$. Note that $\fg_0=\fg^\theta$. 

From \cite[Cor.14]{rlyg} it is known that the following are equivalent: 
\begin{itemize}
\item[(i)] There exists a semisimple element $x\in \fg_1$ for which
$\ad(x):\fg_0\to\fg_1$ is injective;
\item[(ii)] $\theta$ is ell-reg. 
\end{itemize}
Therefore we can also use (i) as the condition for equality in Theorem 1.

Theorem 1 makes no {\it a priori} assumptions on the kinds of elements 
contained in  $\fg_1$. 
But let us now assume that $\fg_1$ contains nonzero semisimple elements. Such gradings are said to have {\bf positive rank}. Their classification is contained in the combined papers \cite{vinberg:graded}, \cite{levy:thetap} and \cite{rlyg}. 

In the case of positive rank gradings, Theorem 1 complements results of 
Panyushev. 
 According to \cite[Thm 2.1]{panyushev:theta}, if $x\in \fg_1$ is semisimple, then 
\begin{equation}\label{pan1}
\frac{1}{m}=\frac{\dim[\fg_0,x]}{\dim[\fg,x]}.
\end{equation}
Since $\dim[\fg_0,x]\leq \dim\fg_0$ with equality exactly when (i) holds for $x$,
and since
$\dim[\fg,x]\leq\dim(\fg/\ft)$ with equality exactly when $x$ is a regular element of $\fg$, 
Theorem 1 combines with \eqref{pan1} to interpose $\dim(\fg/\ft)/m$ in the inequality $\dim[\fg_0,x]\leq \dim\fg_0$. That is, we have
\begin{cor} Assume there is a semisimple element $x\in\fg_1$.
Then we have two inequalities
\[ \dim[\fg_0,x]\underset{1}{\leq} \frac{\dim(\fg/\ft)}{m}\underset{2}{\leq} \dim\fg_0.\]
$\underset{1}{\leq}$ is equality if and only if $x$ is regular (semisimple), and 
$\underset{2}{\leq}$ is equality if and only if $\theta$ is ell-reg. 
\end{cor}
Under the additional assumption that $\fg_1$ contains a  regular semisimple element, Panyushev \cite[4.2]{panyushev:theta} also showed that 
\[\dim\fg_0=\frac{\dim[\fg/\ft]}{m}+k_0,\]
where $k_0\geq 0$ is an integer depending only on the orders $m$ and $e$ of $\theta$ in $\Aut(\fg)$ and $\Out(\fg)$. For example, if $e=1$ then $k_0$ is the number of exponents of $\fg$ divisible by $m$. This is a sharper form of Theorem 1 in the case that $\fg_1$ contains a regular semisimple element.

\subsection{Adjoint Swan Conductors} 

In the setting of section \ref{zeroweight}, sending a representation  $V$ to its highest weight $\lam$ is a simple case of the much broader and still partly conjectural local Langlands correspondence (LLC). 
In section \ref{zeroweight} we saw that the inequalities/equalities of Theorem 1 appear on the dual side of this LLC. 

They also appear on the dual side of the LLC for reductive $p$-adic groups, now as measures of ramification. 

We use notation parallel to that of section \ref{zeroweight}.
Let $k$ be a $p$-adic field and let $G$ be the group of $k$-rational points in a connected and simply-connected almost simple $k$-group $\bG$.  

Let $\wh\fg$ be a simple complex Lie algebra whose root system is dual to that of $\bG$. 

The  LLC predicts the existence of a partition 
\[\Irr^2(G)=\bigsqcup_{\vp}\Pi_\vp,\]
of the set $\Irr^2(G)$ of irreducible discrete series representations of $G$ (up to equivalence) into finite sets $\Pi_\vp$, 
where $\vp$ ranges over certain representations 
$\vp:\scw_k\to \Aut(\wh\fg)$ of the Weil group of $k$. (See \cite{gross-reeder} for more background on the LLC.)  

It is of interest to find invariants relating the discrete series representation $\pi$ of $G$ to the representation $\vp$ of $\scw_k$ for which $\pi\in \Pi_{\vp}$. 

One invariant of $\vp$ is 
its {\it adjoint Swan conductor} $\sw(\vp,\fg)$. 
This is an integer depending only on the image $I=\vp(\sci)$ of the inertia subgroup $\sci\subset \scw_k$. We have a 
factorization $I=S\rtimes P$, where $P$ is a $p$-group and $S$ is a cyclic group of order prime to $p$. We have $\sw(\vp,\fg)\geq 0$, with equality if and only if $P$ is trivial. 

Expected properties of the LLC imply certain inequalities for $\sw(\vp,\fg)$. These inequalities have been found to hold unconditionally. For example if $\vp$ is totally ramified (that is, if $\fg^{I}=0$) then the LLC predicts that
\begin{equation}\label{adswan}
\dim\fg^\theta\leq \sw(\vp,\fg),
\end{equation}
where $\theta$ is a generator of $S$. 
This inequality has been proved in \cite{reeder:adswan} and \cite{bushnell-henniart:swan}. 

Assume now that $p$ does not divide the order of $W$.
Let $m$ be the order of $\theta$. Combining \eqref{adswan} with 
Theorem 1 gives the inequality 
\begin{equation}\label{adswanm}
\frac{\dim\fg/\ft)}{m}\leq \sw(\vp,\fg),
\end{equation}
which is weaker than \eqref{adswan}, but which depends only of the order $m$ of $S$, not on $S$  itself. Moreover, the two inequalities \eqref{adswan} and \eqref{adswanm} coincide if and only if $\theta$ is ell-reg.

\section{Proof of Theorem 1} \label{proof}

The torsion automorphisms of $\fg$ are classified by Kac diagrams. 
We start with a summary of Kac diagrams that I hope is sufficient for the reader to follow the computations. For more background, see \cite{kac:bluebook} and  \cite{reeder:torsion}.

 \subsection{Kac Diagrams}

Fix a divisor $e\in\{1,2,3\}$ of the order of the component group 
$\Out(\fg)$ of $\Aut(\fg)$. Let $\Aut(\fg,e)$ be the set of elements in $\Aut(\fg)$ whose image in $\Out(\fg)$ has order $e$. Then $\Aut(\fg,e)$ has one or two connected components, the latter only when $\fg=\fso_8$ and $e=3$.  

For any torsion automorphism $\theta\in\Aut(\fg,e)$, 
the rank of the fixed point subalgebra $\fg^\theta$ depends only on $e$; we denote this rank by $n_e$. If $e=1$ then $G_1:=\Aut(\fg,1)$ is the identity component of $\Aut(\fg)$ and $n_1$ is is the rank of $\fg$. 

To the pair $(\fg,e)$ one associates an  affine Dynkin diagram $\scd(\fg,e)$.
As we vary over all pairs $(\fg,e)$, the diagrams $\scd(\fg,e)$ range exactly over the affine Coxeter diagrams together with all possible orientations on the multiple edges. If $e=1$, then $\scd(\fg,1)$ is the usual affine Dynkin diagram of $\fg$. 

The vertices in $\scd(\fg,e)$ are indexed by a set $I$ whose cardinality is $n_e+1$, and these vertices are labelled by certain positive integers $\{c_i:\ i\in I\}$, where $1\leq c_i\leq 6$.

The automorphism group $\Aut(\scd(\fg,e))$ of the oriented and labelled diagram $\scd(\fg,e)$ contains a (very small) subgroup $\Om$ with the following property: If $e>1$ then $\Om=\Aut(\scd(\fg,e))$. 
If $e=1$ then $\Om\simeq \pi_1(G)$.

We fix a connected component $\Ga$ of $\Aut(\fg,e)$. For any positive integer $m$, let $\Ga_m$ be the set of elements of $\Ga$ having order $m$. Then $\Ga_m$ is nonempty if and only if $e$ divides $m$. 
The $G_1$-conjugacy classes in $\Ga_m$ are parametrized as follows. 
Let $S_m$ be the set of $I$-tuples $s=(s_i:\ i\in I)$ consisting of   integers $s_i\geq 0$ such that $\gcd\{s_i: i\in I\}=1$ and 
\[m=e\cdot\sum\limits_{i\in I}c_is_i.\]
There is a surjective mapping from $S_m$ to the set of $G_1$-congugacy classes in $\Ga_m$ (Kac coordinates). Two elements $s$ and $ s'\in S_m$ map to the same conjugacy class in $\Ga_m $ if and only if $s$ and $s'$ are conjugate under the group $\Om$. 

For example in $\Ga$ there is a unique conjugacy class of automorphisms of minimal order having abelian fixed-point subalgebras. Such automorphisms are called {\bf principal}. They are ell-reg and have Kac coordinates $s=(s_i)$ where
$s_i=1$ for all $i$. The order of a principal automorphism in $\Ga$, namely  
\[h_e:=e\cdot\sum_{i\in I} c_i,\]
is the Coxeter number of $\Aut(\fg,e)$. 
It is known \cite{reeder:torsion} that equality holds in Theorem 1 for principal elements, namely we have
\begin{equation}\label{coxeternumber}
\frac{1}{h_e}=\frac{n_e}{[\fg:\ft]}.
\end{equation}

For any subset  $J\subset  I$ we set 
\[c_J=\sum_{j\in J} c_j, \qquad c^J=\sum_{i\notin J} c_i.\]
If  $J\neq I$ then the subgraph of $\scd(\fg,e)$ supported on $J$ 
is the finite Dynkin graph of a reductive subalgebra $\fg_J$ of $\fg$. 
Let $|R_J|$ be the number of roots of $\fg_J$.

Let $\theta\in\Ga$ be a torsion automorphism with Kac-coordinates $s=(s_i)$ and let  $J=\{j\in I: s_j=0\}$. 
Then $J\neq I$ and we have $\fg^\theta=\fg_J$. 

\begin{lemma}\label{lemmaJ} The inequality in Theorem 1 for all torsion elements in a component $\Ga\subset\Aut(\fg,e)$ is equivalent to the inequality  
\begin{equation}\label{ineqJ}
n_e\cdot c_J\leq c^J \cdot |R_J|
\end{equation}
for every subset $J\subset I, J\neq I$. 
\end{lemma}
\proof
Let $\theta\in \Ga_m$ have Kac coordinates $(s_i)$ and let $J=\{j\in I: s_j=0\}$. Then
$m\geq e\cdot c^J$ with equality if and only if $s_i=1$ for all $i\in I-J$.
Since 
\[\dim \fg^\theta=\dim \fg_J=n_e+|R_J|, \quad\text{and}\quad
\dim(\fg/\ft)=h_e n_e=e\cdot c_I\cdot n_e
\]
It follows that 
\[\frac{1}{m}\leq \frac{1}{e\cdot c^J}
\quad\text{and}\quad 
\frac{\dim \fg^\theta}{\dim(\fg/\ft)}=\frac{n_e+|R_J|}{e\cdot c_I\cdot n_e}, 
\]
so the inequality in Theorem 1 for every $\theta$  is equivalent to having
\[e\cdot c_I\cdot n_e\leq (n_e+|R_J|)\cdot e\cdot c^J\]
for every $J$. 
Since $c_I= c^J+c_J$, the result follows.
\qed

If $J$ is empty then both sides of \eqref{ineqJ} are zero. We may assume from now on that $J$ is nonempty and that $s_i=1$ for all $i\in I-J$. Thus $J$ is identified with a Kac diagram.

We will show that the integer $f(\fg,e,J)$ defined by 
\[f(\fg,e,J)=c^J|R_J|-n_e c_J 
\]
satisfies $f(\fg,e,J)\geq 0$.  Our analysis will also find  those $J$ for which $f(\fg,e,J)=0$. 

On the other hand, the Kac diagrams of the ell-reg automorphisms of $\fg$ were tabulated in \cite{rlyg} section 7 and are recalled in the Appendix \ref{appendix} of this paper. 
It turns out that the Kac diagrams of ell-reg automorphisms are exactly those for which $f(\fg,e,J)=0$. 

\subsection{Type $A_n$}
The case $\fg=\fsl_{n+1}$ and $e=1$ is very simple but different from the other cases, so we treat it separately here. Fix a nonempty subset $J\subset I$.
 The root system $R_J$ has type 
\[\prod_{i=1}^a A_{q_i}\]
for some positive integers $q_1,\dots,q_a$. Let $q=\sum q_i$. 
Since all $c_i=1$ we have $c_J=q$ and $c^J=n+1-q\geq a$.  Now
\[\begin{split}
f(\fg,1,J)=c^J\sum_{i=1}^a q_i(q_i+1)-(c^J+q-1)q
=c^J\sum_{i=1}^a q_i^2-q^2+q
\geq a\sum_{i=1}^a q_i^2-q^2+q
\geq q.
\end{split}
\]
where the arithmetic-geometric inequality is used in the last step. Since $J\neq \varnothing$ we have
$f(\fg,1,J)\geq q>0$.

\subsection{Preliminary reductions}

In this section, $(\fg,e)$ is of classical type not equal to $(\fsl_{n},1)$. 
We will write $n=n_e$ and $h=h_e$.

The relevant diagrams $\scd(\fg,e)$ for $n\geq 3$ are listed below.  Each diagram has $n+1$ nodes. They are grouped according to their underlying Coxeter diagram.

\[\begin{array}{ccc}
(\fg,e)&\scd(\fg,e)& h=e\cdot c_I\\
\hline
&&\\
\underset{n\geq 2}{ {^2\!A_{2n}}}\quad &\overset{1}\circ\!\!\Longrightarrow\!\!\overset{2}\circ \text{----}\!\! \overset{2}\circ\!\text{--} \cdots\text{--} 
 \overset{2}\circ\!\!\Longrightarrow\!\!\overset{2}\circ
& \quad 2(2n+1)\\
&&\\
 \underset{n\geq 2}{C_n}& \overset{1}\circ\!\!\Longrightarrow\!\!
\overset{2}\circ \text{----}\!\! \overset{2}\circ\!\text{--} \cdots\text{--} 
 \overset{2}\circ\!\!\Longleftarrow\!\!\overset{1}\circ
 &\quad 2n\\
 &&\\
 \underset {n\geq 2}{{^2\!D_{n+1}}}& \overset{1}\circ\!\!\Longleftarrow\!\!
\overset{1}\circ \text{----}\!\! \overset{1}\circ\!\text{--} \cdots\text{--} 
 \overset{1}\circ\!\!\Longrightarrow\!\!\overset{1}\circ
 &\quad 2(n+1)\\
&&\\
\hline
&&\\
\underset {n\geq 3}{{^2\!A_{2n-1}}}\quad&
\begin{matrix}
 \overset{1}\circ \text{----}\!\!\!\!\! &\overset{2}\circ &
 \!\!\!\!\!\text{--} \cdots\text{--}
 \overset{2}\circ\!\!\Longleftarrow\!\!\overset{1}\circ\\
&\underset{1}{\text{\rotatebox{270}{\!\!\!\!\!\!\!\!----\!$\circ$}}}&
 \end{matrix}
 & \quad 2(2n-1) \\
 &&\\
 \underset{n\geq 3}{B_n}&
\begin{matrix}
 \overset{1}\circ \text{----}\!\!\!\!\! &\overset{2}\circ&
 \!\!\!\!\!\text{--} \cdots\text{--}
 \overset{2}\circ\!\!\Longrightarrow\!\!\overset{2}\circ\\
&\underset{1}{\text{\rotatebox{270}{\!\!\!\!\!\!\!\!----\!$\circ$}}}&
 \end{matrix}
 & \quad 2n \\
 &&\\ 
 \hline
 &&\\
 \underset{n\geq 4}{D_n}&
\begin{matrix}
 \overset{1}\circ \text{----}\!\!\!\!\! &\overset{2}\circ&
 \!\!\!\!\!\text{--} \cdots\text{--}\!
 \overset{2}\circ\!\text{----}\overset{1}\circ\\
&\underset{1}{\text{\rotatebox{270}{\!\!\!\!\!\!\!\!----\!$\circ$}}}
&\underset{1}{\text{\rotatebox{270}{\!\!\!\!\!\!\!\!----\!$\circ$}}}
 \end{matrix}
 & \quad 2n-2
 \end{array}
 \]

Let $\scx$ be the set of all triples $(\fg,e,J)$, where $(\fg,e)$ is one of the above classical types for $n\geq 3$ and $J$ is a nonempty proper subset of the set $I$ of vertices of $\scd(\fg,e)$. Let $\scx_0=\{(\fg,e,J):\ f(e,\fg,J)=0\}$. We must prove that $f\geq 0$ on $\scx$ and that $\scx_0$ consists precisely of the diagrams  listed in the appendix \ref{appendix} for classical $(\fg,e)$.

By induction on $n$, we may assume that $f(\fg',e,J')\geq 0$ for all $(\fg',e)$ of the same type as $(\fg,e)$, such that  $\scd(\fg',e)$ has at most $n$ nodes indexed by a set $I'$ and  $J'\subset I'$.

{\bf Definition.\ } If $\scy'\subset \scy$ are subsets of $\scx$, 
we say $\scy'$ is a {\bf refinement} of $\scy$ if 
for every $(\fg,e,J)\in \scy-\scy'$, we have either

(i)\ $f(\fg,e,J)>0$ or 

(ii)\  there exists $(\fg',e',J')\in \scy'$ such that 
$f(\fg,e,J)> f(\fg',e',J')$. 

We observe the following.  
\begin{itemize} 
\item[(i)] Refinement is transitive: if $\scy''$ is a refinement of $\scy'$ and $\scy'$ is a refinement of $\scy$ then $\scy''$ is a refinement of $\scy$.
\item[(ii)] If $\scy$ is a refinement of $\scx$ and 
$f\geq 0$ on $\scy$ then 
$f>0$ on $\scx-\scy$ and $\scx_0=\scy_0$. 
\end{itemize}
From (ii), it suffices to find a refinement $\scy$ of $\scx$ such that $f\geq 0$ on $\scy$ and that $\scy_0$ consists precisely of the ell-reg triples listed in the appendix.
 
Say that  a vertex $i\in I$ is {\it interior} if $i$ is adjacent to at least two other vertices in $\scd(\fg,e)$. In every pair of adjacent vertices, at least one vertex is interior. 
The table of diagrams shows that all interior $i$'s have the same value $c$ of $c_i$ ($c=1$ in type $ {^2D_{n+1}}$ and $c=2$ in the other classical diagrams), and $c\geq c_i$ for all $i\in I$.  

\begin{lemma}\label{contraction} Let $\scy$ be the set of $(\fg,e,J)\in \scx$ for which 
no two interior vertices of $I\!-\! J$ are adjacent in $\scd(\fg,e)$.  Then $\scy$ is a refinement of $\scx$. 
\end{lemma}
\proof 
Consider a triple $(\fg,e,J)\in \scx$
and let $i,j\in I\!-\!J$ be adjacent vertices in $\scd(\fg,e)$. 
 
First assume either that $i$ has degree two, or that $j$ is also interior.  Let $k$ be another vertex adjacent to $i$. The possible configurations of $i,j,k$ in the Kac diagram are
\[\cdots\overset{j}{1}\ \text{---}\ \overset{i}{1}\ \text{---}\ \overset{k}{\ast}\cdots
\qquad\quad\cdots\overset{j}{1}\ \text{---}\ \overset{i}{1}\! =\!=\!\overset{k}{\ast}
\qquad\quad\cdots\overset{k}{\ast}\text{---}\ \overset{i}{1}\! =\!=\!\overset{j}{1}\qquad\qquad
\begin{matrix}
 \overset{k}{\ast}\ \text{---}\!\!\!\! &\overset{i}1&
 \!\!\!\!\!\text{---} 
 \overset{j}1\cdots\\
&\underset{}{\text{\rotatebox{270}{\!\!\!\!\!\!---$\blt$}}}
 \end{matrix}\]
 where $\ast,\blt\in\{0,1\}$ are arbitrary. 
 
Removing $i$ and joining $j$ to $k$ with a bond of the appropriate type, we obtain a diagram $\scd(\fg',e)$ of the the same type as $\scd(\fg,e)$.
The vertices of $\scd(\fg',e)$  are indexed by $I'=I-\{i\}$ and we have $J\subset I'$. In this way, the diagram $\scd(\fg,e,J)$ 
contracts by one vertex to  the diagram $\scd(\fg',e, J)$.
The root system $R'_J$ of $\fg'_J$ is isomorphic to $R_J$, we have $\sum_{i'\in I'-J}c_{i'}=c^J-c$ and $c_J$ is unchanged. 
It follows that 
\[\begin{split}
f(\fg,e,J)-f(\fg',e,J)&=[c^J|R_J|-nc_J]-[(c^J-c)|R_J|-(n-1)c_J]\\
&=c|R_J|-c_J\geq 2c|J|-c|J|=c|J|>0.
\end{split}
\]
Since $|I'-J|=|I-J|-1$, repeating this procedure will eventually produce a diagram $\scd(\fg'',e,J)\in\scy$ and we will have $f(\fg,e,J)>f(\fg'',e,J)$.  

Now assume $i$ has degree three and $j$ is one of the two boundary vertices adjacent to $i$. Let $k$ be the interior vertex adjacent to $i$.

From the first case we may assume $k\in J$. 
So $(\fg,e,J)$ has the form
 
  \[
  \begin{split}
J:\quad \overset{j}{1}\ &\  \overset{i}{1}\ \ \overset{k}{0}\ \ \overset{ q\ \text{vertices}} {\overbrace{0\ \ 0\ \cdots \  0}}\ \ 1\ \cdots \\
 &\ s\qquad\qquad\quad\quad \ \  
\end{split} 
\]
where $s\in \{0,1\}$ and $q\geq 0$. 
Switch $s_i$ and $s_k$ to obtain 
\[
  \begin{split}
J':\quad \overset{j}{1}\ &\  \overset{i}{0}\ \ \overset{k}{1}\ \ \overset{ q\ \text{vertices}} {\overbrace{0\ \ 0\ \cdots \  0}}\ \ 1\ \cdots \\
 &\ s\qquad\qquad\quad\quad \ \  
\end{split} 
\]
If $q=0$ then $f(\fg,e,J)=f(\fg,e,J')>0$ from the previous case. 

Assume $q\geq 1$.
Since $c^{J'}=c^J$,  $n'=n$ and $c_{J'}=c_J$,  we find that 
\[f(\fg,e,J)-f(\fg,e,J')=2(q+s-1)c^J\geq 0\]
with equality only if $q=1$ and $s=0$. In this case,
\begin{equation}\label{vexing}
 \begin{split}
J':\quad \overset{j}{1}\ &\  \overset{i}{0}\ \ \overset{k}{1}
\ 0\ \ 1\ \cdots \\
 &\ 0\qquad\qquad\quad\quad \ \  
\end{split} 
\end{equation}
In section \ref{ACD} we will prove directly that for the case \eqref{vexing} we have  $f(\fg,e,J')>0$, so $f(\fg,e,J)>0$. 
\qed

Our next refinement heads toward equilibrium for the interior components of $R_J$. 

Given a diagram $\scd(\fg,e,J)\in \scx$, 
 let $J^\circ$ be the set of interior vertices in $J$ and let  $R_J^\circ$ be the union of those irreducible components of $R_J$ whose bases are contained in $J^\circ$.  Let   $R_1,\ R_2,\ \dots,\ R_a$ be the components of $R_J^\circ$. Each $R_i$ has type $A_{q_i}$ 
for some integers $q_i\geq 1$.
Let 
\[d(J)=\max\{|q_i-q_j|:\ 1\leq i\leq a\}.\]
be the largest difference between the ranks of any two components of $R_J^\circ$.

\begin{lemma}\label{balance} Let $\scy$ be as in Lemma \ref{contraction} and let $\scz$ be the set of $(\fg,e,J)\in\scy$ for which $d(J)\leq 1$.
Then $\scz$ is a refinement of $\scy$. 
\end{lemma}
\proof 
The value of $f(\fg,e,J)$ is unchanged by permuting the components $R_1,\dots, R_a$. If $d(J)\geq 2$, then we may choose such a permutation to arrange that $q_1-q_2\geq 2$  and there are three interior vertices 
$\{i,j,k\}$ such that $j\in R_1$, $i\in I-J$, $k\in R_2$, as shown. 
\[\cdots\  \overset{j}{0}\ \text{---}\  \overset{i}{1}\ \text{---}\ \overset{k}{0}\ \cdots.
\]
Now switch $s_i$ and $s_j$ to obtain 
a diagram $(\fg,e,J')$: 
\[\cdots\ \overset{j}{1}\ \text{---}\  \overset{i}{0}\ \text{---}\ \overset{k}{0}\ \cdots.
\]
Then $n, c_J$ and $c^J$ are unchanged, and one checks that  
\[f(\fg,e,J)-f(\fg,e,J')=2c^J(q_1-q_2-1)>0.\]
Repeating this process, we eventually find a subset $J''\subset I$ with $f(\fg,e,J)>f(\fg,e,J'')$ and $d(J'')\leq 1$. \qed

It suffices to calculate $f(\fg,e,J)$ only for diagrams $(\fg,e,J)$ in the set $\scz$ of Lemma \ref{balance}. Note that $\scz$ consists of those diagrams for which no two interior vertices in $I-J$ are adjacent and whose components of $R_J^\circ$ have at most two types
$A_{q-1}$ and $A_q$, occuring, say, $x$ and $y$ times respectively.

\begin{lemma} For $(\fg,e,J)\in\scz$, the integer $f(\fg,e,J)$ has the form 
\[f(\fg,e,J)=cxy+\al x+\be y+\ga,\]
where $c$ is the common value of $c_i$ on the interior vertices of $I$ and $\al,\be,\ga$ are polynonomial expressions in $q$.  Also, $\be$ is obtained from $\al$ upon replacing $q$ by $q+1$.
\end{lemma}
\proof
Let  $R_{\partial J}$ be the union of the components of $R_J$ not in $R_J^\circ$ and let 
$\partial J$ be the subset of $J$ supporting $R_{\partial J}$. 
Then 
$R_J=R_{\partial J}\sqcup R_J^\circ$ so we have 
\[|R_J|=|R_{\partial J}|+q(q-1)x+q(q+1)y\quad\text{and}\quad
c_J=c_{\partial J}+c(q-1)x+cqy,
\]
where
\[c_{\partial J}=\sum_{j\in\partial J} c_j.\]
Define integers $a$ and $b$ by 
\[c^J=a+c(x+y)\quad\text{and}\quad n=b+qx+(q+1)y.\]
Then 
\[\begin{split}
f(\fg,e,J)&=|R_J|c^J-nc_J\\
&=[|R_{\partial J}|+q(q-1)x+q(q+1)y]\cdot [a+cx+cy]-
[b+qx+(q+1)y]\cdot c_{\partial J}+c(q-1)x+cqy]\\
&=cxy+\al x+\be y+\ga,
\end{split}
\]
where 
\begin{equation}\label{abc}
\begin{split}
\al&=[c|R_{\partial J}|+aq(q-1)]-[bc(q-1)+qc_{\partial J}]\\
\be&=[c|R_{\partial J}|+aq(q+1)]-[bc(q+1)+(q+1)c_{\partial J}]\\
\ga&=a|R_{\partial J}|-bc_{\partial J},
\end{split}
\end{equation}
as claimed. 
\qed

We will show that $\al,\be,\ga$ are $\geq 0$. 
This implies that $f(\fg,e,J)\geq 0$, with equality if and only if $0=xy=\al=\ga=\be$. 
Without loss of generality, we may then assume $y=0$. Theorem 1 will follow by comparing with the tables  of ell-reg automorphisms. 

\subsubsection{Types ${^2A_{2n}}, C_n, {^2D_{n+1}}$}
\label{ACD}
The underlying Coxeter diagram with indexing set $I=[0,n]$ is 
\[0\!=\!=\!1\ \text{---}\ 2\ \text{---}\ \cdots\ \text{---}\ (n-1)\!=\!=\! n\]
The three types differ only in the labels $c_i$ which do not affect $|R_J|$.
  Let $(\fg, e)$ and $(\fg',e')$ be two of ${^2A_{2n}}, C_n, {^2D_{n+1}}$, with corresponding labellings $c_i, c_i'$. 
For each subset $A\subset I$ we set 
\[c_A=\sum_{i\in A}c_i\qquad c'_A=\sum_{i\in A}c_i'.\]
We set $K=I-J$.

We will compare
\[f=f(\fg,e, J)=|R_J|c_K-nc_J \qquad\text{and}\qquad  
f'=f(\fg',e',J)=|R_J|c'_K-nc'_J .\]
Suppose $(\fg,e)={^2A_{2n}}$ and $(\fg',e')=C_n$. If $n\in K$ then $c_K=c'_K+1$ and $c_J=c_J'$, so $f>f'$. If $n\in J$ then $c_K=c_K'$ and $c_J=c_J'+1$ so $f<f'$. 

Suppose $(\fg,e)={^2A_{2n}}$ and $(\fg',e')={^2D_{n+1}}$. If $0\in K$ then $1+c_K=2c'_K$ and $c_J=2c'_J$, so $2f'>f$. If $0\in J$ then $1+c_J=2c'_J$ and $c_K=2c'_K$, so $f>2f'$. 

Suppose $(\fg,e)={C_n}$ and $(\fg',e')={^2D_{n+1}}$. If $\{0,n\}\subset J$ then 
$c_K=2c_K'$ and $c_J=2c_J'-2$ so $f=2f'+2n>f'$. If $0\in J$ and $1\in K$ then $c_K+1=2c_K'$ and $c_J+1=2c_J'$ so $2f'=f+|R_J|-n$. 
Since no two vertices in $K$ are  adjacent, it follows that $|R_J|\geq n+1$, so $2f'>f$. 

This discussion shows that we need only consider the following three cases.
\[\begin{split}
1. (\fg,e)\ &={{^2A}_{2n}}\quad\text{with}\quad 0\in K\  \text{and}\ n\in J\\
2. (\fg,e)\ &={C_n}\qquad\text{with}\quad \{0,n\}\subset K\\
3. (\fg,e)\ &={{^2D}_{n+1}}\quad\text{with}\quad \{0,n\}\subset J
\end{split}
\]
Indeed, if $f(J,\fg,e)\geq 0$ in cases 1-3 then $f(J,\fg,e)\geq 0$ in all cases and we can only have $f(J,\fg,e)=0$ in cases 1-3.

{\bf Case 1:\ }  $(\fg,e)={^2A_{2n}}$ and $R_J=B_r + xA_{q-1} +y A_q, r\geq 1$. 
\begin{align*}
|R_J|&=2r^2+q(q-1)x+q(q+1)y& 
c_K&= 1+2x+2y\\
n&= r+xq+y(q+1)& 
c_J&= 2(q-1)x+2qy+2r\\
\ga&= 0&
\al&= (q-2r)(q-2r-1)
\end{align*}
Thus we have $f(\fg,e,J)\geq 0$ with equality if and only if $q=2r$ or $2r+1$. These cases are the last two rows in the table of appendix \ref{2A2n} for $n\geq 2$.

{\bf Case 2:\ }  $(\fg,e)={C_n}$ and $R_J= xA_{q-1} +y A_q$. 
\begin{align*}
|R_J|&=q(q-1)x+q(q+1)y& 
c_K&= 2x+2y\\
n&= qx+(q+1)y& 
c_J&= 2(q-1)x+2qy\\
\ga&= 0&
\al&= 0
\end{align*}
Thus we have $f(\fg,e,J)\geq 0$ with equality if and only if $xy=0$. 
These are the cases with $k=x$ in the table of  appendix \ref{Cn}.

{\bf Case 3.\ } $(\fg,e)={^2D_{n+1}}$ and 
$R_J=B_p+xA_{q-1} +y A_q + B_r, \quad r,p>0, q>1.$

\begin{align*}
|R_J|&=2p^2+2r^2+q(q-1)x+q(q+1)y& 
c_K&= 1+x+y\\
n&= p+r+qx+(q+1)y& 
c_J&= p+r+(q-1)x+qy\\
\ga&= (p-r)^2&
\al&= (p-r)^2(p+r-q)(p+r-q+1)
\end{align*}
Thus we have $f(\fg,e,J)\geq 0$ with equality if and only if $xy=0$, $p=r$ and $q=2p$ or $q=2p+1$. These are the cases with $k=q\geq 2$ in the appendix \ref{2Dn}.

\subsubsection{Types ${^2A_{2n-1}}, B_n$ }
The underlying Coxeter diagram with indexing set $I=[0,n]$ is 
\[\begin{split}
1\ \text{---}\  & 2\ \text{---}\ 3\ \text{---}\  \cdots  \ (n-1)\!=\!=\!  n\\
&\underset{\textstyle 0}{\!\text{\rotatebox{270}{\!\!\!\!\!\!\!\!----\!}}}\\
\end{split}
\]

The two types differ only in the label $c_n=1$ for ${^2A_{2n-1}}$ and $c_n=2$ for $B_n$. Comparing as in the previous section we may assume $n\in K$ for ${^2A_{2n-1}}$ and $n\in J$ for $B_n$. 

{\bf Case A1:\ } $\{0,1, n\}\subset J$\quad $R_J=D_p+xA_{q-1}+yA_q$, $p\geq 2$. 
\begin{align*}
|R_J|&=2p(p-1)+q(q-1)x+q(q+1)y& 
c_K&= 1+2x+2y\\
n&= p+qx+(q+1)y& 
c_J&= 2(p-1)+2(q-1)x+2qy\\\
\ga&= 0&
\al&= (2p-q)(2p-q-1)
\end{align*}
In this case we have $f(\fg,e,J)\geq 0$ with equality if and only if $xy=0$ and $q=2p$ or $q=2p-1$. These are the cases with $k=p\geq 2$ in appendix \ref{2A2n-1}.

{\bf Case A2:\ } $\{0,n\}\subset K$, $1\in J$,\quad$
R_J=A_p+xA_{q-1}+yA_q$
\begin{align*}
|R_J|&=p(p+1)+q(q-1)x+q(q+1)y& 
c_K&= 2+2x+2y\\
n&= 1+p+qx+(q+1)y& 
c_J&= 2p-1+2(q-1)x+2qy\\
\ga&= p+1&
\al&= 2(p-q+1)^2+q
\end{align*}
In this case we have $f(\fg,e,J)>0$.

{\bf Case A3: \ } $\{0,1, n\}\subset K$, \quad $R_J=xA_{q-1}+yA_q$, $q\geq 2$
\begin{align*}
|R_J|&=q(q-1)x+q(q+1)y& 
c_K&= 1+2x+2y\\
n&= 1+qx+(q+1)y& 
c_J&= 2(q-1)x+2qy\\
\ga&= 0 &
\al&= (q-1)(q-2)
\end{align*}
In this case we have $f(\fg,e,J)\geq 0$ with equality if and only if $q=2$. 
This is the case $k=n$ in appendix \ref{2A2n-1}.

{\bf Case B1:\ }  $\{0,1\}\subset J$,  \quad $R_J=D_p+xA_{q-1}+yA_q+B_r$ 
\begin{align*}
|R_J|&=2p(p-1)+2r^2+q(q-1)x+q(q+1)y& 
c_K&= 2(1+x+y)\\
n&= p+r+qx+(q+1)y& 
c_J&= 2(p+r-1)+2(q-1)x+2qy\\
\ga&= 2(p-r)(p-r-1) &
\al&= 2(p-r)(p-r-1)+2(p+r-q)^2
\end{align*}
In this case we have $f(\fg,e,J)\geq 0$ with equality if and only if $
p=r$ and $q=2r$, or $p=r+1$ and $q=2r+1$. these are the cases in the last two rows of appendix \ref{Bn} with $k=q$.

{\bf Case B2:\ }  $0\in K,\ 1\in J$, \quad $R_J=A_p+xA_{q-1}+yA_q+B_r$, $p>0$
\begin{align*}
|R_J|&=p(p+1)+2r^2+q(q-1)x+q(q+1)y & 
c_K&=3+2x+2y \\
n&=p+r+1+qx+(q+1)y & 
c_J&=2p+2r-1+2(q-1)x+2qy \\
\ga&=(2r-p-1)^2+3r & 
\al&=2(p-q)(p-q+1)+(q-r)^2+3r^2+2p
\end{align*}
In this case we have $f(\fg,e,J)>0$.

{\bf Case B3:\ } $\{0,1\}\subset K$, \quad  $R_J=xA_{q-1}+yA_q+B_r$, $r\geq 1$.
\begin{align*}
|R_J|&=2r^2+q(q-1)x+q(q+1)y & 
c_K&=2+2x+2y \\
n&=r+1+qx+(q+1)y & 
c_J&=2r+2(q-1)x+2qy \\
\ga&=2r(r-1) & 
\al&=2(q-r-1)^2+2r(r-1)
\end{align*}
In this case we have $f(\fg,e,J)\geq 0$ with equality if and only if $r=1$ and $q=2$. This is the case $k=2$ in appendix \ref{Bn}

\subsubsection{Type $D_n$}

We choose the indexing set $I=\{0,1,\dots, n\}$ as in \cite{bour456}, so that  $I_{1}=\{0,1,n-1,n\}$. 
Up to automorphisms of $\scd$, there are five cases for $J\cap I_{1}$. 

{\bf Case 1: } $\{0,1,n-1,n\}\subset J$ \quad $
R_J=D_p\times x A_{q-1}\times y A_q\times D_r,\quad p,q,r\geq 2$. 
\begin{align*}
|R_J|&=2p(p-1)+2r(r-1)+q(q-1)x+q(q+1)y & 
c_K&=2+2x+2y \\
n&=p+r+qx+(q+1)y & 
c_J&=2(p+r-2+(q-1)x+qy) \\
\ga&=2(p-r)^2& 
\al&=2(p-r)^2+2(p-q+r)(p-q+r-1)
\end{align*}

In this case we have $f(\fg,e,J)\geq 0$ with 
equality holds if and only if $p=r$ and $q=2p$ or $q=2p-1$. These are  the cases $2<k=q$ in appendix \ref{Dn}

{\bf Case 2: } $\{0,n\}\subset J, \{1,n-1\}\subset K$, 
$R_J=A_{p-1}+xA_{q-1}+yA_q+A_{r-1}$, $p,r\geq 2$. 

\begin{align*}
|R_J|&=p(p-1)+r(r-1)+q(q-1)x+q(q+1)y & 
c_K&=2(2+x+y) \\
n&=p+r+qx+(q+1)y & 
c_J&=2(p+r-3+(q-1)x+qy) \\
\ga&=2(p-r)^2+2(p+r)& 
\al&=2(p-q)^2+2(p-r)^2+2q
\end{align*}
In this case $f(\fg,e,J)>0$.

{\bf Case 3: } $\{0,1\}\subset J, \{n-1,n\}\subset K$, 
$R_J=D_p+xA_{q-1}+yA_q$, $p\geq 2$.

\begin{align*}
|R_J|&=2p(p-1)+q(q-1)x+q(q+1)y & 
c_K&=2(1+x+y) \\
n&=1+p+qx+(q+1)y & 
c_J&=2(p-1+(q-1)x+qy) \\
\ga&=2(p-1)^2& 
\al&=2[ (p-q+1)^2+(p-2)(p-1) + (q-2)]
\end{align*}
In this case $f(\fg,e,J)>0$. 

{\bf  Case 4: } $0\in J, \{1, n-1,n\}\subset K$, $R_J=A_{p-1}+xA_{q-1}+yA_q$, 
\begin{align*}
|R_J|&=p(p-1)+q(q-1)x+q(q+1)y & 
c_K&=3+2x+2y) \\
n&=1+p+qx+(q+1)y & 
c_J&=2p-3+2(q-1)x+2qy \\
\ga&=(p-1)^2+2& 
\al&=2(p-q)^2+(q-1)^2+1]
\end{align*}
In this case $f(\fg,e,J)>0$. 

{\bf Case 5: } $\{0,1, n-1,n\}\subset K$, $R_J=xA_{q-1}+yA_q$, $q\geq 2$.
\begin{align*}
|R_J|&=q(q-1)x+q(q+1)y & 
c_K&=2+2x+2y) \\
n&=2+qx+(q+1)y & 
c_J&=2(q-1)x+2qy \\
\ga&=0& 
\al&=2(q-1)(q-2)]
\end{align*}
In this case $f(\fg,e,J)\geq 0$ with equality if and only if $q=2$. This is the case $k=2$ in appendix \ref{Dn}.

\section{Exceptional Lie Algebras} \label{exceptional}

In principal (or on a computer) one can verify Theorem 1 for the exceptional Lie algebras and ${^3D_4}$ by checking the theorem for each subset $J\subset I$. The aim of this section is to make this verification somewhat more illuminating.  

Assume the diagram $\scd(\fg,e)$, with labels $c_i$ has one of the types 
\[\begin{array}{cc}
{G_2}\quad &
\overset{1}\circ\text{----}\overset{2}\circ\!\Rrightarrow\!\overset{3}\circ \\
&\\
{^3\!D_{4}}\quad &
\overset{1}\circ\text{----}\overset{2}\circ\!\Lleftarrow\!\overset{1}\circ \\
&\\
{F_4} \quad&
\overset{1}\circ\!\text{----}
\overset{2}\circ\!\text{----}
\overset{3}\circ\! \Longrightarrow\!\!
\overset{4}\circ\text{----}
\overset{2}\circ \\
&\\
{^2\!E_{6}}\quad &
\overset{1}\circ\!\text{----}\!
\overset{2}\circ\! \text{----}\!
\overset{3}\circ\!\!\Longleftarrow\!\!
\overset{2}\circ\text{----}\!
\overset{1}\circ\\
\end{array}
\]
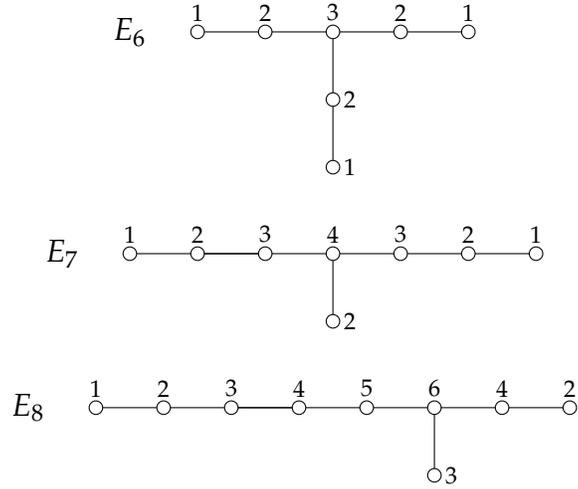
\begin{center}
\begin{tikzpicture}[scale=.9]
\node at (-3,1) {$E_6$};
\draw (0,-1) circle [radius=.1] node[right]{$\scriptstyle 1$}; 
\draw (0,0) circle [radius=.1] node[right]{$\scriptstyle 2$};
\draw (0,1) circle [radius=.1] node[above]{$\scriptstyle 3$};
\draw (-2,1) circle [radius=.1] node[above]{$\scriptstyle 1$};
\draw (-1,1) circle [radius=.1] node[above]{$\scriptstyle 2$};
\draw (1,1) circle [radius=.1] node[above]{$\scriptstyle 2$};
\draw (2,1) circle [radius=.1] node[above]{$\scriptstyle 1$};
\draw (0,-.9) -- (0,-.1);
\draw (0,.1) -- (0,.9);
\draw (.1,1) -- (.9,1);
\draw (1.1,1) -- (1.9,1);
\draw (-.1,1) -- (-.9,1);
\draw (-1.1,1) -- (-1.9,1);
\end{tikzpicture}

\begin{tikzpicture}[scale=.9]
\node at (-4,1) {$E_7$};
\draw (-3,1) circle [radius=.1] node[above]{$\scriptstyle 1$};
\draw (-2.9,1) -- (-2.1,1);
\draw (-2,1) circle [radius=.1] node[above]{$\scriptstyle 2$};
\draw (-1.9,1) -- (-1.1,1);
\draw (-1,1) circle [radius=.1] node[above]{$\scriptstyle 3$};
\draw (-.9,1) -- (-.1,1);
\draw (.1,1) -- (.9,1);
\draw (0,1) circle [radius=.1] node[above]{$\scriptstyle 4$};
\draw (1.1,1) -- (1.9,1);
\draw (1,1) circle [radius=.1] node[above]{$\scriptstyle 3$};
\draw (2.1,1) -- (2.9,1);
\draw (2,1) circle [radius=.1] node[above]{$\scriptstyle 2$};
\draw (-1.9,1) -- (-1.1,1);
\draw (3,1) circle [radius=.1] node[above]{$\scriptstyle 1$};
\draw (0,0) circle [radius=.1] node[right]{$\scriptstyle 2$};
\draw (0,.9) -- (0,.1);
\end{tikzpicture}

\begin{tikzpicture}[scale=.9]
\node at (-5,1) {$E_8$};
\draw (-4,1) circle [radius=.1] node[above]{$\scriptstyle 1$};
\draw (-3.9,1) -- (-3.1,1);
\draw (-3,1) circle [radius=.1] node[above]{$\scriptstyle 2$};
\draw (-2.9,1) -- (-2.1,1);
\draw (-2,1) circle [radius=.1] node[above]{$\scriptstyle 3$};
\draw (-1.9,1) -- (-1.1,1);
\draw (-1,1) circle [radius=.1] node[above]{$\scriptstyle 4$};
\draw (-.9,1) -- (-.1,1);
\draw (.1,1) -- (.9,1);
\draw (0,1) circle [radius=.1] node[above]{$\scriptstyle 5$};
\draw (1.1,1) -- (1.9,1);
\draw (1,1) circle [radius=.1] node[above]{$\scriptstyle 6$};
\draw (2.1,1) -- (2.9,1);
\draw (2,1) circle [radius=.1] node[above]{$\scriptstyle 4$};
\draw (-1.9,1) -- (-1.1,1);
\draw (3,1) circle [radius=.1] node[above]{$\scriptstyle 2$};
\draw (1,0) circle [radius=.1] node[right]{$\scriptstyle 3$};
\draw (1,.9) -- (1,.1);
\end{tikzpicture}
\end{center}

\subsection{Small $J$}
We begin with cases where $|J|\leq 8$. 

When $R_J=A_1$, Theorem 1 follows from an observation which applies uniformly to all exceptional cases. 
Namely, each coefficient $c_i$ is at most twice the average of the remaining coefficients, with equality just for the largest coefficient $c_i=c$ whose node is the target of the arrow or is the branch node. On the other hand, the  Kac diagrams
\[11\Rrightarrow 0\quad 10\Lleftarrow 1\quad 
111\Rightarrow01\quad 110\Leftarrow 11\quad
\EVI{1}{1}{0}{1}{1}{1}{1}\qquad 
\EVII{1}{1}{1}{0}{1}{1}{1}{1}\qquad
\E{1}{1}{1}{1}{1}{0}{1}{1}{1}
\]
are those of the ell-reg automorphisms of order $h-ec$. 

Now suppose $R_{J}=2A_1$. Then $J=\{i,\ j\}$ where $i,j$ are not adjacent in $\scd(\fg,e)$. The maximum value of $c_i+c_j$ is $2c-2$ (with $c$ as above), which gives 
\[|R_J|c^J-nc_J\geq 2(n-2c+4)\geq 0
\]
with equality only in $G_2, F_4, E_8$. On the other hand, the Kac diagrams
\[01\Rrightarrow 0\qquad 
101\Rightarrow01\qquad 
\E{1}{1}{1}{0}{1}{0}{1}{1}{1}\]
are those of ell-reg automorphisms of order $h-2c+2$.

If $R_J=A_2$, one finds similarly that  
\[|R_J|c^J-nc_J=6c^I-(n+6)(c_i+c_j)\geq 0,
\]
with equality only in ${^3D_4}$. The Kac diagram 
\[00\Lleftarrow 1\]
Is the ell-reg automorphism of order $e=3$.

If $R_J=B_2$ or $G_2$, one finds that $|R_J|c^J-nc_J>0$. 

At this point, the theorem is proved for $G_2$ and ${^3D_4}$, and we may assume 
$R_J$ has rank at least three in the remaining cases. 

Assume that $R_J=3A_1$. Then $f(\fg,e,J)=6c^I-(n+6)c_J$. The Kac diagrams with maximal $c_J$ are 
\[ 010\Rightarrow10\qquad 010\Leftarrow 10\qquad
\EVI{1}{0}{1}{0}{1}{0}{1}\qquad 
\EVII{1}{0}{1}{0}{1}{0}{1}{1}\qquad
\E{1}{0}{1}{0}{1}{0}{1}{1}{1}.
\]
These all have $f(\fg,e,J)\geq 0$, with equality only in the $E_6$ case. This is the Kac diagram of the ell-reg inner automorphism of $\fg=E_6$ of order six.

Assume that $R_J=A_1+A_2$. In the same manner we find $f(\fg,e,J)\geq 0$, 
with equality only in the cases
\[ 101\Rightarrow00,\quad 100\Leftarrow 10,
\]
which are the Kac diagrams for the ell-reg automorphisms of $F_4$ of order four and the outer ell-reg automorphism of $E_6$ of order six. 

Assume that $R_J=4A_1$. This only exists in type $E$ 
we find $f(\fg,e,J)\geq 0$, 
with equality only in the case
\[\E{0}{0}{1}{0}{1}{0}{1}{}{1}.\]
This is the ell-reg automorphism of $E_8$ of order 15.

\subsection{$F_4$ and ${^2E_6}$}
We now complete the proof of Theorem 1 for $(\fg,e)$ of type $F_4$ and ${^2E_6}$, for which $\scd(\fg,e)$ has the same underlying Coxeter diagram. By the previous section, we may assume $|R_J|>8$. 
Arguing as in section \ref{ACD} we need only consider cases of the form
\[ \ast\ \ast\ \ast\Rightarrow 00\qquad 
\ast\ \ast\ \ast\Leftarrow 11 \qquad  
\ast\ \ast\ \ast\Leftarrow 10\qquad \ast\ \ast\ \ast\Leftarrow 01,
\]
with $R_J>8$. The possibilities are

\[\begin{array}{c|c|l}
J& R_J\cdot c^J& 4\cdot c_J\\
\hline
100\Rightarrow 00&48\cdot 1 &4\cdot 11 \\
010\Rightarrow 00&20 \cdot 2&4\cdot10 \leftarrow \\
001\Rightarrow 00&12 \cdot 3&4\cdot9 \leftarrow \\
110\Rightarrow 00&18\cdot 3 &4\cdot9 \\
\hline
000\Leftarrow 11 &  12\cdot 3 & 4\cdot 6\\
\hline
000\Leftarrow 10 &  14\cdot 2 & 4\cdot 7\leftarrow\\
\hline
000\Leftarrow 01 &  32\cdot 1 & 4\cdot 8\leftarrow\\
100\Leftarrow 01 &  18\cdot 2 & 4\cdot 7\\
010\Leftarrow 01 &  10\cdot 3 & 4\cdot 6
\end{array}
\]
We have $f(\fg,e,J)\geq 0$ with equality in the cases marked by $\leftarrow$. 
These are the elliptic regular automorphisms of order $2$ for $F_4$ and outer automorphisms of $E_6$ of orders $4$ and $2$. This completes the proof of Theorem 1 in the cases $F_4$ and ${^2E_6}$.

\subsection{Types $E_6,\ E_7,\ E_8$}

We consider the ends of the interval $1<m<h$ in two steps. 

{\bf Step 1.\ } For  each $1<m<n$ we compute the minimum
\[r(m)=\min\{|R_J|:\ c^J=m\}.\]
In the tables below we check that 
\begin{equation}\label{step1}
 r(m)\geq \frac{|R|}{m}-n
 \end{equation}
for each $m<n$ and we verify that equality holds in \eqref{step1} 
for at most one $J$ with 
$c^J=m$. This will prove Theorem 1 when $m<n$.

Next we consider $|R_J|$ where $c^J\geq n$. If $|R_J|> h-n$, then since $h=|R|/n$ and $c^J\geq n$ we have 
\[|R_J|> h-n=\frac{|R|}{n}-n\geq \frac{|R|}{c^J}-n.
\]
Hence we may also assume $|R_J|\leq h-n$. Since we have already proved Theorem 1 for $|R_J|\leq 8$, we may in fact assume 
\[10\leq r<h-n\]

{\bf Step 2.\ } For  each even integer $r\leq h-n$ we compute the minimum 
\[m(r)=\min\{c^J:\ |R_J|=r\}.\]
In the tables below we check that 
\begin{equation}\label{step2}
 r\geq \frac{|R|}{m(r)}-n
 \end{equation}
for each $r\leq h-n$ and we verify that equality holds in \eqref{step2} for at most one $J$ with
$|R_J|=r$. This will complete the proof of Theorem 1.  

\subsubsection{$E_6$}
Step 1 for $E_6$ is tabulated as follows. The bold entries are those types of $R_J$ for which 
$r(m)=|R_J|$. 
\[
\begin{array}{c|l|c|c|c}
m& \text{types of  $R_J$ with $c^J=m$}&r(m)& (|R|/m)-6&\text{$J$}\\
\hline
2&{\bf A_1A_5},\ \ D_5&32 &30&\text{none}\\
\hline
3&{\bf 3A_2},\ \ A_1A_4,\ D_4,\ A_5&18&18
&\EVI{0}{0}{1}{0}{0}{0}{0}\\
\hline
4&{\bf 2A_2A_1},\ {\bf A_1A_3},\ \ 2A_1A_3,\ A_4&14&12
&\text{none}\\
\hline
5&{\bf 2A_1A_2},\ 2A_2A_1,\ \ A_1A_3,\ A_3&10&42/5
&\text{none}
\end{array}
\]
Since $h-n=12-6<8$, the proof for $E_6$ is completed by Step 1 alone. 

\subsubsection{$E_7$}

Step 1 for $E_7$ is tabulated as follows. 
\[\begin{array}{c|l|c|c|c}
m& \text{types of  $R_J$ with $c^J=m$}&r(m)& (|R|/m)-7&\text{$J$}\\
\hline
2&{\bf A_7},\ \ A_1D_6,\ \ E_6&56 &56&\EVII{0}{0}{0}{0}{0}{0}{0}{1}\\
\hline
3&{\bf A_2A_5},\ \ A_1D_5,\ \ A_6,\ \ D_6,\ &36&35
&\text{none}\\
\hline
4&{\bf 2A_3A_1},\ \ {\bf A_2A_4},\ \ 2A_1D_4,\ \ A_5,\ \ A_1A_5,\ \ D_5 &26&49/2
&\text{none}\\
\hline
5&{\bf A_1A_2A_3},\ \ A_1A_4,\ \ A_2A_4, \ \ A_1D_4, \ \ A_5,\ \ A_1A_5 &20&91/5
&\text{none}\\
\hline
6&{\bf 2A_2A_1},\ 2A_1A_3,\ A_2A_3,\ 3A_2,\  3A_1A_3,&14&14
&\EVII{1}{0}{0}{1}{0}{0}{1}{0}\\
& \ A_4,\ A_1A_4,\ 2A_3,\ D_4,\ A_5&&
\end{array}
\]
For step 2, we need only consider $r=10$. The only root systems with 10 roots are  $5A_1$ and $2A_1A_2$. All occurrences of these as $R_J$ in $E_7$ have $c^J=8$. Since  
$|R|/8-7<10$, Theorem 1 is now proved for $E_7$. 

\subsubsection{$E_8$}

In Step 1, we take $1< m< 8$ and compute $r(m)$ in the following table.
The types of $R_J$ for which $m_J=m$ are shown; those for which $|R_J|= r(m)$ are in bold face. The right column gives the unique $J$ for which  $r(m_J)=(240/m)-8$, if it exists. 
\[\begin{array}{c|l|c|c|c}
m& \text{types of}\  R_J:\ m_J=m& r(m)& (240/m)-8&J\\
\hline
2&{\bf D_8},\ \ A_1E_7&112&112&
\E{0}{\ 0 \ }{\  0\ }{\ 0\ }{\  0\ }{\ 0}{\ \ 0\ }{\ 1\ }{ 0}\\
\hline
3&{\bf A_8},\ A_2E_6,\ D_7,\ \ E_7&72&72&
\E{0}{\ 0 \ }{\  0\ }{\ 0\ }{\  0\ }{\ 0}{\ \ 0\ }{\ 0\ }{ 1}\\
\hline
4&{\bf A_3D_5},\ A_7,\ A_1A_7,\ A_1D_6,\ \ A_1A_6&52&52&
\E{0}{\ 0 \ }{\  0\ }{\ 1\ }{\  0\ }{\ 0}{\ \ 0\ }{\ 0\ }{ 0}\\
\hline
5&{\bf 2A_4},\ A_1A_6,\ A_2D_5,\  A_7,\ D_6,\ A_1E_6 &40&40&
\E{0}{\ 0 \ }{\  0\ }{\ 0\ }{\  1\ }{\ 0}{\ \ 0\ }{\ 0\ }{ 0}\\
\hline
6&{\bf A_3A_4},\ A_1^2A_5,\  A_3D_4,\  A_2A_5,\ A_1A_2A_5,&32&32&
\E{1}{\ 0 \ }{\  0\ }{\ 0\ }{\  1\ }{\ 0}{\ \ 0\ }{\ 0\ }{ 0}\\
&A_1D_5,\ \ A_6,\  \  A_1^2D_5 ,\  A_7,\ \ E_6&&&\\
\hline
7&{\bf A_1A_2A_4},\ A_2D_4,\ A_3A_4, \ A_1A_5&28&184/7&\text{none} \\
&A_1D_5,\ A_6,\  A_1A_6, \ A_2D_5  &&& \\
\end{array}
\]

For Step 2, we take $r=10,12,\dots, 22$ and compute $m(r)$ in the following table.
The types of $R_J$ for which $|R_J|=r$ are shown; those for which $c^J=m(r)$ are in bold face and that $J$ for which $|R_J|=(240/m_J)-n$, if it exists, is shown in the right column. 
\[\renewcommand{\arraystretch}{1.2}
\begin{array}{c|l|c|c|c}
r& \text{types of}\  R_J \ \text{with}\ |R_J|=r& m(r)& [240/m(r)]-8&J\\
\hline
10&{\bf A_1^5},\  {\bf A_1^2A_2}&14&64/7&\text{none}\\
\hline
12&{\bf A_1^3A_2},\  A_2^2,\  A_3&12&12&
\E{1}{\ 0 \ }{\  1\ }{\ 0\ }{\  0\ }{\ 1}{\ \ 0\ }{\ 1\ }{ 0}\\
\hline
14&{\bf A_1^4A_2},\  A_1A_2^2,\  A_1A_3&12&12&\text{none}\\
\hline
16&{\bf A_1^2A_2^2},\  A_1^2A_3&10&16&
\E{1}{\ 0 \ }{\  1\ }{\ 0\ }{\  0\ }{\ 1}{\ \ 0\ }{\ 0\ }{ 0}\\
\hline
18&{\bf A_2A_3},\  A_1^3A_3,\ A_2^3  &10&16&\text{none}\\
\hline
20&{\bf A_1A_2A_3},\ {\bf A_1A_2^3},\ A_4\  &9&56/3&\text{none}\\
\hline
22&{\bf A_1^2A_2A_3},\   A_1A_4&8&22&
\E{0}{\ 1 \ }{\  0\ }{\ 0\ }{\  0\ }{\ 1}{\ \ 0\ }{\ 0\ }{ 0}
\end{array}
\]
In each case we have $r\geq[240/m(r)]-8$ and equality is 
achieved by at most one $J$, as indicated in the rightmost column.

The proof of Theorem Theorem 1 for $E_8$ is now complete.

\newpage
\section{Appendix: The classification of ell-reg automorphisms }\label{appendix}

For reference in the proofs above, we recall the classification of ell-reg automorphisms given in \cite{rlyg}. There is only one inner ell-reg automorphism of $\fsl_n$, namely the principle one, so we ignore this case. 

\subsection{Type ${^2A_{2n}}$}\label{2A2n}

The ell-reg outer automorphisms of $\fsl_{2n+1}$ correspond to odd quotients $d$ of $2n$ and $2n+1$. The graphs $\scd(\fsl_{2n+1}, 2)$ are as shown:

\[\overset{n=1}{\twoAtwo}\qquad\qquad 
\overset{n>1}{
 \overset{1}\circ\!\!\Longrightarrow\!\!\overset{2}\circ \text{----}\!\! \overset{2}\circ\!\text{--} \cdots\text{--} \overset{2}\circ\!\text{----}
 \overset{2}\circ\!\!\Longrightarrow\!\!\overset{2}\circ}
\]
The ell-reg outer automorphisms of $\fsl_{2n+1}$ correspond to odd quotients $d$ of $2n$ and $2n+1$. 
We write these quotients as 
\[d=\frac{2n+1}{2k+1}\qquad d=\frac{n}{k}\]
respectively. The cases overlap only when $d=1$.
The corresponding ell-reg automorphism has order $m=2d$ in both cases. 
\begin{center}
\[
{\renewcommand{\arraystretch}{1.3}
\begin{array}{ c c }
\hline
d=m/2&s\\
\hline
3&\twoAtwoo\\
%\hline
2 &\twoAtwooo \\
\hline
 \end{array}}
 \]
 \end{center}
\begin{center}
\[
{\renewcommand{\arraystretch}{1.7}
\begin{array}{ cc c }
\hline
d=m/2&\quad&s \\
\hline
2n+1
&&1 \Rightarrow 1\ \ 1\ \cdots\ 1\ \ 1\ \Rightarrow 1 \\
%\hline
1&&1 \Rightarrow 0\ \ 0\ \cdots\ 0\ \ 0\ \Rightarrow 0 \\
%\hline
\dfrac{2n+1}{2k+1} &&
1 \Rightarrow\underset{A_{2{k}} } { \underbrace{ 0\cdots0}  }\ \ 1\ \ 
\underset{A_{2{k}} } { \underbrace{ 0\cdots0}  }\ \ 1\ \ \cdots\ \ 1\ \ 
\underset{B_k}  { \underbrace{ 0\cdots0  \Rightarrow 0}}
\\
%\hline
\dfrac{n}{k}\quad \quad &&
1 \Rightarrow\underset{A_{2{k}-1} } { \underbrace{ 0\cdots0}  }\ \ 1\ \ 
\underset{A_{2{k}-1} } { \underbrace{ 0\cdots0}  }\ \ 1\ \ \cdots\ \ 1\ \ 
\underset{B_{k} } { \underbrace{ 0\cdots0 \Rightarrow 0}}
\\
\hline
 \end{array}}
 \]
 \end{center}
In the two last rows we have $0<k<n$ such that $d$ is odd and the number of type-$A$ factors is $(d-1)/2$. The next-to-last row corrects an error in \cite{rlyg}.

\subsection{Type ${^2A_{2n-1}}$}\label{2A2n-1}
The graph $\scd(\fsl_{2n},2)$, with  $n\geq 3$ and labels $c_0,\ c_1,\ \dots\ ,\ c_n$ is shown here,with $c_0=c_n=1$. 
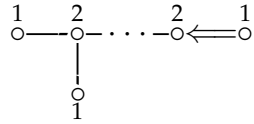
\[
\begin{matrix}
 \overset{1}\circ \text{----}\!\!\!\!\! &\overset{2}\circ&
 \!\!\!\!\!\text{--} \cdots\text{--}
 \overset{2}\circ\!\!\Longleftarrow\!\!\overset{1}\circ\\
&\underset{1}{\text{\rotatebox{270}{\!\!\!\!\!\!\!----$\circ$\!}}}&
 \end{matrix}
\]

The ell-reg outer automorphisms of $\fsl_{2n}$  correspond to odd quotients $d$ of $2n-1$ and $2n$. 
We write these quotients as 
\[d=\frac{2n-1}{2k-1}\qquad d=\frac{n}{k}\]
respectively. The cases overlap only when $d=1$.
The corresponding ell-reg automorphism has order $m=2d$ in both cases. 

\begin{center}
\[
{\renewcommand{\arraystretch}{1.2}
\begin{array}{ cc c}
\hline
d=m/2&\quad &s\\
\hline
&&\\
2n-1&&
\begin{split}
1\ \  \ & 1\ \ \ 1\ \ \ 1\ \cdots 1\ \ \ 1\Leftarrow1\\
&1
\end{split}\\
&&\\
1&&
\begin{split}
0\ \  \ & 0\ \ \ 0\ \ \  0\cdots 0\ \ \ 0\Leftarrow1\\
&0
\end{split}
\\
&&\\
\underset{\text{$n$ odd}}{n}&&
\begin{split}
1\ \  \ & 0\ \ \ 1\ \ \ 0\ \ \ 1\cdots 1\ \ \ 0\Leftarrow1\\
&1\\
\end{split}\\
&&\\
\dfrac{2n-1}{2k-1}&&
\begin{split}
&
\overset{ D_{{k}}     } {\overbrace{0\ \ \ \ 0\ \cdots \ 0}}\ \ 1\ \ 
\overset{A_{2k-2} } {\overbrace{ 0\ \cdots \ 0}}\ \ 1\cdots
\ 1\ \ \overset{A_{2k-2} } {\overbrace{ 0\ \cdots\ 0}}
\Leftarrow 1\\
&\quad \ \ 0\\
\end{split} 
\\
&&\\
\dfrac{n}{k}&&
\begin{split}
&\overset{ D_{{k}}     } {\overbrace{0\ \ \ \ 0\ \cdots \ 0}}\ \ 1\ \ 
\overset{A_{2{k}-1} } {\overbrace{ 0\ \cdots \ 0}}\ \ 1\cdots
\ 1\ \ \overset{A_{2{k}-1} } {\overbrace{ 0\ \cdots\ 0}}
\Leftarrow 1\\
&\quad \ \ 0\\
\end{split} 
\\
\hline
\end{array}}
\]
\end{center}
In the last two rows we have $1<k<n$ such that $d$ is odd and there are $(d-1)/2$ components of type $A$.

\newpage
\subsection{Type $B_n$}\label{Bn}

The graph $\scd(\fso_{2n+1},1)$ with labels $c_0,\ c_1,\ \dots\ ,\ c_n$ is shown here,with $c_0=c_n=1$. 
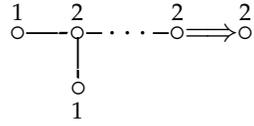
\[
\begin{matrix}
 \overset{1}\circ \text{----}\!\!\!\!\! &\overset{2}\circ&
 \!\!\!\!\!\text{--} \cdots\text{--}
 \overset{2}\circ\!\!\Longrightarrow\!\!\overset{2}\circ\\
&\underset{1}{\text{\rotatebox{270}{\!\!\!\!\!\!\!\!----$\circ$}}}&
 \end{matrix}
\]
The ell-reg automorphisms of $\fso_{2n+1}$  are of the form 
$\pi^k$, where $\pi$ is a principal automorphism and $k$ is a divisor of $n$. The order $m$ of $\pi^k$ is $m=2n/k$ and the Kac coordinates of $\pi^k$ are given in the table below. We replace each node $i$ by the Kac coordinate $s_i\in\{0,1\}$, 
and also omit the single bonds in the graph. Recall that $J=\{i\in I:\ s_i=0\}$. 
\begin{center}
\[
{\renewcommand{\arraystretch}{1.3}
\begin{array}{ccc}
k\mid n&m &s=(s_0,\ s_1,\ \cdots, s_n)\\
\hline 
1&2n&
\begin{split}
1\ \ \ & 1\ \ \ 1\ \ \ 1\ \ \ 1\cdots 1\ \ \ 1\Rightarrow1\\
&1
\end{split}\\
\hline
2 &n&
\begin{split}
1\ \ \ & 0\ \ \ 1\ \ \ 0\ \ \ 1\cdots 0\ \ \ 1\Rightarrow0\\
&1
\end{split}
\\
\hline
\underset{ \text{$k$ even} }{k>2} &\dfrac{2n}{k}&
\begin{split}
&\overset{ D_{{k}/2}     } {\overbrace{0\ \ \ 0\ \cdots \ 0}}\ \ \ 1\ \ \ 
\overset{A_{{k}-1} } {\overbrace{ 0\ \cdots \ 0}}\ \ 1\ \ 0\cdots0\ \ 
1\ \ \ \overset{A_{{k}-1} } {\overbrace{ 0\ \cdots\   0}}\ \ \ 1\ \ \ 
\overset{B_{{k}/2} } {\overbrace{ 0\cdots 0\  \Rightarrow 0}}\\
&\quad \ 0
\end{split} 
\\
\hline
\underset{\text{$k$ odd} }{k>1} &\dfrac{2n}{k}&
\begin{split}
&\overset{ D_{{(k+1)}/2}     } {\overbrace{0\ \ \ 0\ \cdots \ 0}}\ \ \ 1\ \ \ 
\overset{A_{{k}-1} } {\overbrace{ 0\ \cdots \ 0}}\ \ 1\ \ 0\cdots0\ \ 
1\ \ \ \overset{A_{{k}-1} } {\overbrace{ 0\ \cdots\   0}}\ \ \ 1\ \ \ 
\overset{B_{{(k-1)}/2} } {\overbrace{ 0\cdots 0\  \Rightarrow 0}}\\
&\quad \ 0
\end{split}
\\
\end{array}}
\]
\end{center}
The second line $m=n$ only occurs if $n$ is even.
In the last two lines there are $\dfrac{n}{k}-1$ factors of type $A_{k-1}$. 
\subsection{Type $C_n$}\label{Cn}
The graph $\scd(\fsp_{2n},1)$ with labels $c_0,\ c_1,\ \dots\ ,\ c_n$ is shown here,
with $c_0=c_n=1$. 
The Coxeter number is $2n$.
\[
 \overset{1}\circ\!\!\Longrightarrow\!\!\overset{2}\circ \text{----}\!\! \overset{2}\circ\!\text{--} \cdots\text{--} 
 \overset{2}\circ\!\!\Longleftarrow\!\!\overset{1}\circ
\]
As with $\fso_{2n+1}$, the ell-reg automorphisms of $\fsp_{2n}$ are powers $\pi^k$ of a principle automorphism $\pi$, where $k$ is a divisor of $n$.  The order $m$ of $\pi^k$ is $m=2n/k$ and the Kac coordinates of $\pi^k$ are given in the table below.
\begin{center}
\[
{\renewcommand{\arraystretch}{1.3}
\begin{array}{ccc}
k\mid n &m& s=(s_0,\ s_1,\ \cdots, s_n)\\
\hline 
1&2n&
1\Rightarrow 1\ \ \ 1\ \ \ 1\ \ \ 1\cdots 1\ \ \ 1\Leftarrow1
\\
\hline
k>1 &\frac{2n}{k}
&
1 \Rightarrow\underset{A_{{k}-1} } { \underbrace{ 0\cdots0}  }\ \ 1\ \ 
\underset{A_{{k}-1} } { \underbrace{ 0\cdots0}  }\ \ 1\ \ \cdots\ \ 1\ \ 
\underset{A_{{k}-1} } { \underbrace{ 0\cdots0}  }\Leftarrow 1 
\\

\end{array}}
\]
\end{center}
In the last line for $k>1$ there are $n/k$ factors of type $A_{k-1}$,

\subsection{Type $D_n$} \label{Dn}

The graph $\scd(\fso_{2n},1)$ with labels $c_0,\ c_1,\ \dots\ ,\ c_n$ is shown here,
with $c_0=c_1=c_{n-1}=c_n=1$.

\begin{center}

\begin{tikzpicture}[scale=.9]
%\node at (-3,1) {$E_6$};
\draw (-3,0) circle [radius=.1] node[above]{$\scriptstyle 1$}; 
\draw (-2,0) circle [radius=.1] node[above]{$\scriptstyle 2$}; 
\draw (-1,0) circle [radius=.1] node[above]{$\scriptstyle 2$}; 
\draw(-.9,0) --(-.5,0);
\node at (0,0) {$\cdots$};
\draw(.5,0) --(.9,0);
\draw (1,0) circle [radius=.1] node[above]{$\scriptstyle 2$};
\draw (2,0) circle [radius=.1] node[above]{$\scriptstyle 2$};
\draw (3,0) circle [radius=.1] node[above]{$\scriptstyle 1$};
\draw (-2,-1) circle [radius=.1] node[below]{$\scriptstyle 1$};
\draw (2,-1) circle [radius=.1] node[below]{$\scriptstyle 1$};
\draw (-2.9,0) -- (-2.1,0);
\draw (-1.9,0) -- (-1.1,0);
\draw (1.1,0) -- (1.9,0);
\draw (2.1,0) -- (2.9,0);
\draw (-2,-.1) -- (-2,-.9);
\draw (2,-.1) -- (2,-.9);
\end{tikzpicture}
\end{center}

The ell-reg congugacy classes in $\Aut(\fso_{2n},1)$ correspond to even divisors $k$ of $n$ (where $m=2n/k$) and odd divisors $k$ of $n-1$ (where $m=(2n-2)/k$), as shown in the table below.
\begin{center}
\[ 
{\renewcommand{\arraystretch}{1.3}
\begin{array}{ccc}
k&m& s=(s_i)\\
\hline 
1&2n-2\quad&
\begin{split}
1\ \ \ \!& \ 1\ \ \ 1\cdots\  1\ \ \  1\ \ \ 1\\
&\ 1\qquad\qquad\ \ \!1
\end{split}\\
\hline
2&\underset{\text{$n$ even}}{n}&
\begin{split}
1\ \ & 0\ \ 1\ \ 0\ \ 1\cdots1\ \  0\ \ 1\ \  0\ \ 1\\
&1\qquad\qquad\quad\qquad\ \ \ \ \!\!1
\end{split}\\
\hline
n& \underset{\text{$n$ even}}{2}&
 \begin{split}
&\overset{ D_{n/2}   } {\overbrace{0\ \ 0\cdots \  0}}\ \ 1\ \ 
\overset{D_{n/2} } {\overbrace{ 0\ \cdots\  0\  \  \ 0}}\\
 &\quad 0\qquad\qquad\quad\quad \ \  0\\
\end{split} \\
\hline
 \textstyle
    \begin{array}{c}
    \text{$k$ even}\\
    \text{$k$ divides $n$}\\
      2<k<n
    \end{array} &\dfrac{2n}{k}\quad 
   &
 \begin{split}
&\overset{ D_{{k}/2}     } {\overbrace{0\ \ 0\cdots \  0}}\ \ 1\ \ 
\overset{A_{{k}-1} } {\overbrace{ 0\  \cdots \  0}}\ \ 1\ 0\ \cdots\ 
0\ 1\ \ \overset{A_{{k}-1} } {\overbrace{ 0\ \cdots\  0}}\ \ 1\ \ 
\overset{D_{{k}/2} } {\overbrace{ 0\ \cdots\  0\  \  \ 0}}\\
&\quad 0\qquad\qquad\qquad\qquad\qquad\qquad\qquad\qquad\qquad\quad
\ \ \ \ \!0\\
\end{split} 
\\
\hline
\textstyle
    \begin{array}{c}
    \text{$k$ odd}\\ \text{$k$ divides\ $ n-1$}\\
      1<k<n-1
    \end{array}& \dfrac{2n-2}{k}\quad 
&
 \begin{split}
&\overset{ D_{({k}+1)/2}     } {\overbrace{0\ \ 0\cdots \ 0}}\ \ 1\ \ 
\overset{A_{{k}-1} } {\overbrace{ 0\  \cdots \ 0}}\ \ 1\ 0\ \cdots
\ 0\ 1\ \ \overset{A_{{k}-1} } {\overbrace{ 0\  \cdots\   0}}\ \ 1\ \ 
\overset{D_{({k}+1)/2} } {\overbrace{ 0\ \cdots\  0\  \  \ 0}}\\
&\quad 0\qquad\qquad\qquad\qquad\qquad\qquad\qquad\qquad\qquad\quad\ \ \ \ \! 0
\end{split} 
\\
\end{array}}
\]
\end{center}

In the last two rows, the number of type -$A$ factors is one less than $\frac{n}{k}$ and ${n-1}{k}$, respectively.

\subsection{Type ${^2D_{n+1}}$}\label{2Dn}
The graph $\scd(\fso_{2n+2},2)$ (with $n\geq 2$) with  is shown here,
with $c_0=c_1=\cdots=c_n=1$. 
\[{^2D_{n+1}}:\quad 
 \overset{1}\circ\!\!\Longleftarrow\!\!\overset{1}\circ \text{----}\!\! \overset{1}\circ\!\text{--} \cdots\text{--}\!  \overset{1}\circ\!\! \text{----}
 \overset{1}\circ\!\!\Longrightarrow\!\!\overset{1}\circ
\]

The ell-reg classes in $\Aut(\fso_{2n+2},2)$ correspond to even divisors $k$ of $n$ with order $m=2n/k$ and odd divisors $k$ of $n+1$ with order $2(n+1)/k$.

\begin{center}
\[
{\renewcommand{\arraystretch}{1.5}
\begin{array}{ccc}
\hline
k&m&s=(s_0,s_1,\dots, s_n)\\
\hline
1&2n+2&
1\Leftarrow 1\ \ 1\cdots1\ \ \ 1 \Rightarrow 1
\\
&&\\
\hline
2&
\underset{\text{$n$ even}}{n}&
0\Leftarrow 1\ \ 0\ \ 1\ \ 0\cdots 0\ \ \ 1\ \ 0\ \  1 \Rightarrow 0
\\
&&\\
\hline
\textstyle
    \begin{array}{c}
    \text{$k$ even}\\ \text{$k$ divides\ $ n$}\\
      2<k
    \end{array}& \dfrac{2n}{k}\quad 
&
\overset{ B_{{k}/2}     } {\overbrace{0\Leftarrow 0\cdots  0}}
\ \ 1\ \ 
\overset{A_{{k}-1} } {\overbrace{ 0\cdots  0}}\ \ 1\cdots
\ 1\  \overset{A_{{k}-1} } {\overbrace{ 0 \cdots 0}}\ \ 1\ \ 
\overset{B_{k/2} } {\overbrace{ 0\cdots 0 \Rightarrow 0}}
\\
\hline
\textstyle
    \begin{array}{c}
    \text{$k$ odd}\\ \text{$k$ divides\ $ n+1$}\\
      1<k
    \end{array}& \dfrac{2n+2}{k}\quad&
\overset{ B_{(k-1)/2}} {\overbrace{0\Leftarrow 0\cdots  0}}
\ \ 1\ \ 
\overset{A_{{k}-1} } {\overbrace{ 0\cdots 0}}\ \ 1\cdots
\ 1\ \overset{A_{{k}-1} } {\overbrace{ 0\cdots 0}}\ \ 1\ \ 
\overset{B_{({k}-1)/2} } {\overbrace{ 0\cdots 0 \Rightarrow 0}}\\
\hline
\end{array}}
\]
\end{center}
In the last two rows, the number of type -$A$ factors is one less than $\frac{n}{k}$ and $\frac{n+1}{k}$, respectively.

\subsection{Exceptional Lie Algebras}\

\[
\begin{array}[t]{cc}
&E_6\\
\hline
m&s\\
\hline 
&\\
12&\EVI{1}{1}{1}{1}{1}{1}{1}\\
&\\
 9&\EVI{1}{1}{1}{1}{0}{1}{1} \\
&\\
6&\EVI{1}{ 1}{ 0}{ 0}{ 1}{ 0}{ 1}\\
 &\\
3&\EVI{0}{ 0}{ 0}{ 0}{ 1}{ 0}{ 0}\\
\hline
\end{array}\qquad
\begin{array}[t]{cc}
&{^2E_6}\\
\hline
m&s\\
\hline
&\\
18&\outEVI{1}{1}{1}{1}{1}\\
 12&\outEVI{1}{1}{0}{1}{1}\\
 6&\outEVI{1}{0}{0}{1}{0} \\
4&\outEVI{0}{0}{0}{1}{0}\\
2&\outEVI{0}{0}{0}{0}{1} \\
\hline
\end{array}\qquad
\begin{array}[t]{cc}
&E_7\\
\hline
m&s\\
\hline
&\\
18&\EVII{1}{1}{1}{1}{1}{1}{1}{1}\\
 14&\EVII{1}{1}{1}{1}{0}{1}{1}{1}\\
 6&\EVII{1}{ 0}{ 0}{ 0}{ 1}{ 0}{ 0}{ 1}\\
2&\EVII{0}{ 0}{ 1}{ 0}{ 0}{ 0}{ 0}{ 0} \\
\hline
\end{array}
\qquad
\begin{array}[t]{cc}
&E_8\\
\hline
m&s\\
\hline
&\\
30&\E{1}{ 1}{ 1}{ 1}{ 1}{ 1}{ 1}{ 1}{ 1} \\
 24&\E{1}{ 1}{ 1}{ 1}{ 0}{ 1}{ 1}{ 1}{ 1}\\
20&\E{1}{ 1}{ 1}{ 1}{ 0}{ 1}{ 0}{ 1}{ 1}\\
15&\E{1}{ 1}{ 0}{ 0}{ 1}{ 0}{ 1}{ 0}{ 1}\\
12&\E{1}{ 1}{ 0}{ 0}{ 1}{ 0}{ 0}{ 1}{ 0}\\
10&\E{1}{ 0}{ 0}{ 0}{ 1}{ 0}{ 0}{ 1}{ 0} \\
8&\E{0}{ 0}{ 0}{ 0}{ 1}{ 0}{ 0}{ 0}{ 1}\\
6&\E{1}{ 0}{ 0}{ 0}{ 0}{ 1}{ 0}{ 0}{ 0} \\
5&\E{0}{ 0}{ 0}{ 0}{ 0}{ 1}{ 0}{ 0}{ 0}\\
4&\E{0}{ 0}{ 0}{ 0}{ 0}{ 0}{ 1}{ 0}{ 0} \\
3&\E{0}{ 0}{ 1}{ 0}{ 0}{ 0}{ 0}{ 0}{ 0}\\
2&\E{0}{ 1}{ 0}{ 0}{ 0}{ 0}{ 0}{ 0}{0} \\
\hline
\end{array}
\]

\[
\begin{array}[t]{cc}
&F_4\\
\hline
m&s \\
\hline
12&\FIV{1}{1}{1}{1}{1}\\
8&\FIV{1}{1}{1}{0}{1}\\
6&\FIV{1}{0}{1}{0}{1}\\
4&\FIV{1}{0}{1}{0}{0}\\
3&\FIV{0}{0}{1}{0}{0}\\
2&\FIV{0}{1}{0}{0}{0}\\
\hline
\end{array}
\qquad\qquad\qquad\quad
\begin{array}[t]{cc}
&G_2\\
\hline
m&s \\
\hline
6&1\ 1\Rrightarrow 1 \\
3&1\ 1\Rrightarrow 0\\
2&0\ 1\Rrightarrow 0\\
\hline
\end{array}
\qquad\qquad\qquad\quad
\begin{array}[t]{cc}
&{^3D_4}\\
\hline
m&s \\
\hline
12&1\ 1\Lleftarrow 1 \\
6&1\ 0\Lleftarrow1\\
3&1\ 0\Lleftarrow 0\\
\hline
\end{array}
\]

\def\noopsort#1{}
\providecommand{\bysame}{\leavevmode\hbox to3em{\hrulefill}\thinspace}


\begin{thebibliography}{10}

\bibitem{ariki}
S. ~Ariki, 
\emph{On the representations of the Weyl groups of type $D$ corresponding to the zero weights of the representations of $\SO(2n,\bc)$.}
Proc. Sym. Pure Math. \textbf{47} (1987) no. 2,  pp.~327--342 

\bibitem{aik}
I. ~Amemiya, N. ~ Iwahori, K. ~Koiki,
\emph{On some generalization of B. Kostant's partition function},
Manifolds and Lie groups, Papers in Honor of Yoz\^o Matsushima, Prog. Math., \textbf{14}, Birkh\"auser (1981)


\bibitem{amt}
S. ~Ariki, J. ~ Matsuzawa, I. ~Terada,
\emph{Representations of Weyl groups on Zero Weight spaces of $\fg$-Modules},
Algebraic and Topological Theories - to the memory of Dr. Takehiko Miyata pp.~546--568 (1985)

%\bibitem{atiyah-bott}
%\emph{A Lefschetz Fixed Point Formula for Elliptic Complexes: II. Applications}, 
%Annals Math., \textbf{88} (3) 1968 pp.~451--491.

\bibitem{atiyah-segal}
M.~Atiyah, G.~Segal,
\emph{The Index of Elliptic Operators II}, 
Annals Math., \textbf{87} (3) pp.~531--545.


\bibitem{beck-robins} 
M.~Beck, S.~Robins
\emph{Computing the continuous discretely},  
Springer-Verlag (2007).

%\bibitem{borel:131} 
%A.~Borel, 
%\emph{Properties and linear representations of Chevalley groups},  
%Seminar in  Algebraic Groups and related finite groups, Lecture Notes in Math. \textbf{131}, Springer-Verlag (1970), 
 % pp.~1--55.






%\bibitem{borel:linear}
%A.~Borel,
%\emph{Linear algebraic groups},
%Graduate Texts in Mathematics,
%vol.~126,
%Springer-Verlag New York, 1991.

%\bibitem{borel-serre:sousgroupes}
%A.~Borel, J.P.~Serre,
%\emph{Sur certains sous-groupes des groupes de Lie compacts},
%Comment. Math. Helv., \textbf{27} (1953), pp.~128-139.

%\bibitem{bour123}
%N.~Bourbaki,
%\emph{Lie groups and Lie algebras}, Chap. 1-3, Springer-Verlag, Berlin, 2002.


\bibitem{bour456}
N.~Bourbaki,
\emph{Lie groups and Lie algebras}, Chap. 4-6, Springer-Verlag, Berlin, 2002.

\bibitem{bour76}
N.~Bourbaki,
\emph{Lie groups and Lie algebras}, Chap. 7-9, Springer-Verlag, Berlin, 2002.


%\bibitem{brion-vergne:partition}
%M.~Brion and M.~Vergne,
%\emph{Residue formulae, vector partition functions and lattice points in rational polytopes}, 
%J. Amer. Math. Soc., \textbf{10} (4) 1997 pp.~797--833.

% \bibitem{bruhat:boulder}
%F.~Bruhat, 
%\emph{$\fp$-adic groups},  Proc. Symp. Pure  Math., \textbf{9},    Amer. Math. Soc., Providence, RI (1966), 
%  pp.~63--70. 




%\bibitem{bruhat-tits}
%F.~Bruhat and J.~Tits, \emph{Groupes r\'eductifs sur un corps 
%local, I, II}, 
%Publ. Math. Inst. Hautes \'Etudes Sci. 
%No. 48, (1972 ), pp.~5--251, 
%No. 60, (1984), pp. ~197--376.

%\bibitem{bushnell-henniart:GL(2)}
%C.~Bushnell, G.~Henniart
%\emph{The Local Langlands conjecture for $GL(2)$}, Springer-Verlag, Berlin, 2006.

\bibitem{bushnell-henniart:swan}
C.~Bushnell, G.~Henniart
\emph{Tame multiplicity and conductor for local Galois representations}, Tunisian J. Math. 2(2) (2020), pp.~337-357.

\bibitem{ecartan:1927}
\'E.~Cartan
\emph{La g\'eom\'etrie des groupes simples},
Ann. di Mat., \textbf{4} (1927), pp.~209-256.


%\bibitem{carter:finite}
%R.~Carter, 
%\emph{Finite groups of Lie type. Conjugacy 
%classes and complex characters}, Wiley, 1985.

%\bibitem{carter:weyl} 
%R.~Carter,  
%\emph{Conjugacy classes in the Weyl group},  
%Compositio Math., \textbf{25}  (1972), 
 % pp.~1--59.
 
%\bibitem{cartier:weyl} 
%P.~Cartier,  
%\emph{On H. Weyl's character formula},  
%Bull. Amer. Math. Soc., \textbf{67}, no.2  (1961), 
% pp.~228--230.


%\bibitem{debacker-reeder:Lpackets}
%S.~DeBacker and M.~Reeder, 
%\emph{Depth-Zero Supercuspidal $L$-packets and their Stability},
%Annals of. Math., \textbf{169} no.3  (2009), pp.~795--901.



\bibitem{gross-reeder}
B. ~Gross, M. ~Reeder,
\emph{Arithmetic invariants of discrete Langlands parameters},
Duke Math. Jour., \textbf{154} (2010), pp.~431--508.

\bibitem{kac:bluebook}
V.~Kac,
\emph{Infinite dimensional Lie algebras}, 3rd ed.,
Cambridge, 1995.

%\bibitem{kac:legendre}
%\bysame,
%\emph{Simple Lie groups and the Legendre Symbol}, Lecture Notes in Math., \textbf{848}, 
%(1981) pp.~110-124.

%\bibitem{kazhdan-lusztig:affinefixedpoints}
%D.~Kazhdan and G.~Lusztig,
%\emph{Fixed point varieties on affine flag manifolds}, 
%Israel Jour. Math., \textbf{62} (1988), pp.~129-168.

\bibitem{kostant:betti}
B. ~Kostant, 
\emph{The principal three-dimensional subgroup and the 
Betti numbers of a complex simple Lie group},
Amer. J. Math., \textbf{81} (1959), pp.~973--1032.

%\bibitem{kostant:weight}
%\bysame
%\emph{A formula for the multiplicity of a weight},
%Trans. AMS., \textbf{93} (1959), pp.~53--73.

%\bibitem{kostant:poly}
%\bysame
%\emph{Lie Group Representations on Polynomial Rings},
%Amer. J. Math., \textbf{85}, No. 3 (1963), pp.~327--404.

\bibitem{kostant:eta}
\bysame
\emph{On Macdonald's  $\eta$-function formula, the Laplacian and Generalized Exponents},
 Adv. Math., \textbf{20} (1976),  
pp.~179--212.

\bibitem{levy:thetap}
\bysame,
\emph{Vinberg's $\theta$-groups in positive characteristic and Kostant-Weierstrass slices},
Transform. Groups, \textbf{14}, no. 2, (2009), pp.~417--461.

\bibitem{panyushev:theta}
D. ~Panyushev, 
\emph{On invariant theory of $\theta$-groups},
Jour. Algebra, \textbf{283} (2005), pp.~655--670.


\bibitem{prasad:cox}
D.~Prasad,
\emph{Half the sum of positive roots, the Coxeter element, and a theorem of Kostant},
arxiv:1402.5504

%\bibitem{reeder:cohom}
%M.~Reeder, 
%\emph{On the cohomology of compact Lie groups},
%L'Enseignement Math., \textbf{41} (1995), %pp.~181--200.


%\bibitem{reeder:prehom}
%\bysame,
%\emph{Whittaker functions, prehomogeneous vector spaces and standard representations of $p$-adic groups},
%J. Reine angew. Math., \textbf{520 } (1997), pp.~37--93.

%\bibitem{reeder:elliptic}
%\bysame, 
%\emph{Euler Poincare pairings and elliptic representation theory of Weyl groups and $p$-adic groups},
%Comp. Math., \textbf{129} (2001), pp.~149--181.



%\bibitem{reeder:cyclotomic}
%\bysame, 
%\emph{Elliptic centralizers in Weyl groups and their coinvariant representations},
%Representation Theory, \textbf{15} (2011), pp.~63--111.

%\bibitem{reeder:leveltwo}
%\bysame,
%\emph{Level-two structure of simply-laced Coxeter groups},
%Jour. Alg., \textbf{285}(2005), pp.~29--57.

\bibitem{reeder:adswan}
\bysame,
\emph{Adjoint Swan conductors I: The essentially tame case},
Int. Math. Research Notices, \textbf{2018} (9) (2018), pp.~2661--2692.


\bibitem{reeder:torsion}
\bysame,
\emph{Torsion automorphisms of simple Lie algebras},
L'Enseignement Math., \textbf{56}(2) (2010), pp.~3--47.

\bibitem{reeder:weyl}
\bysame,
\emph{Weyl group characters afforded by zero weight spaces},
preprint.



\bibitem{rlyg}
M. ~Reeder, P. Levy,  J.-K. Yu, B. ~Gross,
\emph{Gradings of positive rank on simple Lie algebras},
Transformation Groups, \textbf{17}, No 4,  (2012), pp.~1123--1190.

  


\bibitem{springer:regular}
T.A.~Springer, \ 
\emph{Regular elements in finite reflection groups},  Inv. Math.,
\textbf{25} (1974), pp.~159--198.


% \bibitem{sturmfels:partition}
% B.~Sturmfels,\ 
%\emph{On Vector Partition Functions},
%J. Comb. Thy., \textbf{72} (1995), pp.~302--309. 
   
%\bibitem{varadarajan}
%V.S. ~Varadarajan, 
%\emph{Lie groups, Lie algebras, and their representations},  
%Springer-Verlag, New York,  1984.
   
%\bibitem{vinberg:greenbook}
%E.B. ~Vinberg, 
%\emph{Linear representations of groups},  Birkhauser  1989.

 
%\bibitem{wigner}
%E. P. ~Wigner, 
%\emph{Group Theory and its application to the quantum mechanics of atomic 
%spectra},  Academic Press, New York,  (1959),
 % pp.~209--224. 
 
\bibitem{vinberg:graded}
E.B.~Vinberg, 
\emph{The Weyl group of a graded Lie algebra}, Izv. Akad. Nauk SSSR Ser. Mat., 
\textbf{40} no. 3 (1976), pp.~488-526. English translation: Math. USSR-Izv. \textbf{10} (1977), 
pp.~463--495.




\end{thebibliography}
\end{document}